\documentclass[12pt]{amsart}
\usepackage{ amsthm,amsfonts,amssymb,graphics,graphicx,fullpage,amsmath,timestamp}
\numberwithin{equation}{section}

\usepackage[latin1]{inputenc}\relax

\newcommand{\ha}{{\hat \alpha}}

\newcommand{\btheta}{\bar \theta}

\newcommand{\Tpsi}{{\tilde{\psi}}}

\newcommand{\Fyx}{{\CalF}^{y\to x}\,}

\newcommand{\CalF}{\mathcal{F}}

\newcommand{\CalS}{\mathcal{S}}
\newcommand{\CalU}{\mathcal{U}}

\newcommand{\Span}{\text{\rm span}}
\newcommand{\trans}{\text{\rm tr}}

\newcommand{\mA}{{\mathbb A}}
\newcommand{\RR}{{\mathbb R}}

\newcommand{\CC}{{\mathbb C}}

\newcommand\cD{${\mathcal  D}$}

\newcommand\br{\begin{rem}}
\newcommand\er{\end{rem}}
\newcommand\bp{\begin{pmatrix}}
\newcommand\ep{\end{pmatrix}}
\newcommand\be{\begin{equation}}
\newcommand\ee{\end{equation}}
\newcommand\ba{\begin{equation}\begin{aligned}}
\newcommand\ea{\end{aligned}\end{equation}}

\newcommand\adots{\mathinner{\mkern2mu\raise1pt\hbox{.}
\mkern3mu\raise4pt\hbox{.}\mkern1mu\raise7pt\hbox{.}}}

\newtheorem{theo}{Theorem}[section]
\newtheorem{prop}[theo]{Proposition}
\newtheorem{cor}[theo]{Corollary}
\newtheorem{lem}[theo]{Lemma}
\newtheorem{defi}[theo]{Definition}

\newtheorem{rem}[theo]{Remark}
\newtheorem{rems}[theo]{Remarks}

\numberwithin{equation}{section}

\title{Pointwise Green Function Bounds\\ and stability of
combustion waves}
\date{Last Updated:  \timestamp}

\author{Gregory Lyng, Mohammadreza Roofi,\\
Benjamin Texier, and Kevin Zumbrun}

\thanks{
K.Z. was partially supported  by NSF grant DMS-0300487.
B.T. was partially supported by NSF grant DMS-0505780.
}

\begin{document}

\begin{abstract}
Generalizing similar results for viscous shock and relaxation waves, we
establish sharp pointwise Green function bounds and linearized
and nonlinear stability for traveling wave solutions of an
abstract viscous combustion model
including both Majda's model and the full reacting compressible
Navier--Stokes equations with artificial viscosity with general
multi-species reaction and reaction-dependent equation of state,
under the necessary conditions of
strong spectral stability, i.e., stable point
spectrum of the linearized operator about the wave, transversality of
the profile as a connection in the traveling-wave ODE,
and hyperbolic stability  of the associated Chapman--Jouguet
(square-wave) approximation.
Notably, our results apply to combustion waves of any type:
weak or strong, detonations or deflagrations, reducing the study
of stability to verification of a readily numerically checkable
Evans function condition.
Together with spectral results of Lyng and Zumbrun, this gives
immediately stability of small-amplitude strong detonations in
the small heat-release (i.e., fluid-dynamical) limit, simplifying
and greatly extending previous results obtained by energy methods
by Liu--Ying and Tesei--Tan for Majda's model and the
reactive Navier--Stokes equations, respectively.
\end{abstract}

\maketitle

\section{Introduction}

In this paper, we extend the viscous shock stability theory of \cite{ZH, MaZ2, MaZ3, MaZ4, HZ} to traveling waves of combustion models, including the simplified combustion model of Majda, and an artificial viscosity version of the reacting Naver-Stokes equations.
Specifically, we (i) derive sharp pointwise Green function
bounds, yielding a sharp $L^1 \cap L^p \to L^p$ linearized stability
criterion in terms of an Evans function condition, and
(ii) assuming the Evans stability condition, establish
nonlinear stability for waves of arbitrary type: weak or strong
detonation, weak of strong deflagration.

The results described in this paper, Theorems \ref{D} and \ref{nonlin} below, represent in particular the first stability results of any
kind for large-amplitude combustion waves and for weak detonations of
Majda's model.

This reduces the question of linear and nonlinear stability to
verification of the Evans condition, a criterion that is readily
checked numerically \cite{Br1, Br2, Br3, BrZ, BDG}.

\begin{subsection}{Combustion models}
We show that viscous shock and combustion waves,
like their hyperbolic counterparts, can be studied within a common framework.
Indeed, viscous shocks, viscous detonations, and relaxation shocks may
all be considered as traveling waves of the special class of
hyperbolic--parabolic balance laws, or reaction--diffusion--convection
equations,
\begin{equation}\label{RDC}
U_t + {\mathcal F } (U) = 0, \qquad {\mathcal F} (U) = F(U)_x-(B(U)U_x)_x - G(U),
\end{equation}
having the damping property
\be\label{damping}
\Re \,\sigma(dG)\le 0,
\ee
where (here and elsewhere)
$\sigma(M)$ denotes spectrum of a matrix or linear operator $M$.
For viscous shocks, $G\equiv 0$, while for relaxation shocks,
$dG$ has constant rank, its kernel corresponding to a local
equilibrium manifold.

By contrast, combustion equations have the composite structure
$$
G(U)=\phi(U)\tilde G(U),
$$
where $\phi$ is a scalar ``ignition function'' that turns the
reaction on or off-- specifically, it is zero on some subset of
the state space and positive elsewhere-- and $\tilde G$ is a
relaxation type term,
 $d\tilde G$ has constant rank
  and
$\sigma(d\tilde G)\le 0$: that is, an interpolation between the
viscous and relaxation case.
Thus, traveling combustion waves exhibit features of {\it both}
viscous and relaxation shocks, in various different regimes, and
our analysis must take this into account.

Specifically, we study a subclass of \eqref{RDC}, comprising systems of the form,
\begin{equation} \label{model} \left\{ \begin{aligned}
u_t+f(u,z)_x& =b u_{xx}+qk\phi(u)z, \\
z_t & = d z_{xx}-k\phi(u)z,
\end{aligned} \right. \end{equation}
where $u \in \RR^n$ and $z \in \RR^r,$ and $\phi$ is a ``bump''-type ignition function.
  The physical constant $q$ is the heat release parameter.
Here, $q>0$ corresponds to an exothermic reaction.

When $n=1$ and $r = 1,$ \eqref{model} is Majda's single-reaction
combustion model. Then, $u$ is a lumped variable combining various
aspects of specific volume, particle velocity, and temperature,
while $z\in[0,1]$ is the mass fraction of reactant. The positive
constant $k$ represents the rate of the reaction. In Majda's
model, the diffusion coefficients $b$ and $d$ are also assumed to
be positive constants. In the following, we scale the variables so
that
$b \equiv 1.$

The vectorial version of \eqref{model}, with
$u\in \RR^n$, $z\in \RR^r$, and $b$ and $d$ positive definite matrices,
is sufficient to encompass the artificial viscosity version
of the full reactive compressible
Navier--Stokes equations written in Lagrangian coordinates,
with multi-species reaction and reaction-dependent equation of state,
where $u=(\tau, v, E)$, with $\tau$, $v$, and $E$ denoting
specific volume, velocity, and energy density, $\phi=\phi(T)$,
$z_1, \dots, z_r$ denoting mass fractions of reactant species,
and $k$ matrix-valued with eigenvalues of strictly negative real part
\cite{ZKochel, LyZ2}.

Throughout the paper,
we shall carry out in parallel the analysis of the scalar and
the (artificial viscosity) system case,
exposing the main ideas in the simpler setting
of Majda's model, then indicating by a series of brief remarks
the extension to the general case.

Physical (as opposed to artifical) diffusion terms are of form
$(b(u,z)u_x)_x$ and $(d(u,z)z_x)_x$ with $b$, $d$ matrix-valued
and $b$ semi-definite \cite{LyZ1, LyZ2}. The diffusion
coefficient $b$ is commonly assumed to depend on $u$ alone;
however, like the equation of state, it properly
depends on the make-up of the gas,
hence on the mass fraction $z$ of the reactant. See comments in section \ref{comments}, about the extension of the results of this paper to such systems.

\end{subsection}


\begin{subsection}{Statement of the results}

Consider a general traveling-wave solution $U(x,t)=\bar U(x-st)$
of \eqref{RDC}, and the associated linearized equation about $\bar U$
in moving coordinates $(x-st,t)$:
\be\label{genlin}
U_t + L U = 0, \qquad L U := ({\mathcal F}'(\bar U) - s \partial_x) U.
\ee

 Proposition \ref{RHexist} below, adapted from \cite{LyZ2}, discusses the question of the existence of traveling waves for the scalar version of \eqref{model};
regarding the system version, see Remark \ref{sysconnection}.2.

\begin{defi} Let $X$ and $Y$ be two Banach spaces. A traveling wave $\bar U$ solution of \emph{\eqref{RDC}} is said to be $X \to Y$ linearly orbitally stable
if, for any solution $\tilde U$ of \emph{\eqref{genlin}} with initial data in $X,$ there exists a phase shift $\delta,$ such that $\tilde U(\cdot,t)$ approaches $\delta(t)\bar U'(\cdot),$ in $Y$ and as $t\to \infty.$
\end{defi}

Our first main theorem is the following linearized stability criterion,
precisely analogous to that of the viscous shock case, in terms of
the Evans function $D(\cdot)$
associated with the linearized operator $L$ about the wave, an analytic
function defined on frequencies $\lambda:\, \Re \lambda\ge 0$,
whose zeroes correspond to eigenvalues of $L$ (see Section \ref{evans}
for further details).

\begin{theo}\label{D}
A traveling combustion wave of \emph{\eqref{model}}
is $L^1\cap L^p\to L^p$ linearly orbitally stable for $p>1$ if and only if
\begin{equation} \label{D-condition}
  D(\cdot) \mbox{  has precisely one zero in $\{\Re \lambda \ge 0\}$ (necessarily at $\lambda=0$)}.
\end{equation}
\end{theo}

Theorem \ref{D} is obtained as a result of detailed pointwise bounds
on the Green function of the linearized equations about the wave,
given in Proposition \ref{greenbounds} (resp. Remark \ref{outgoing}.4
in the system case); see Section \ref{linstab}.

Note that the spatial derivative $\bar U'(x)$ of the traveling-wave profile
$\bar U$ is always a zero eigenfunction of $L$, a consequence of
translation invariance of the original evolution equation.
Thus, at a formal level,
condition \eqref{D-condition} corresponds to the statement that perturbations in all
directions other than translation decay with time to the order
of linear approximation, or linearized orbital stability.
(At the formal level only, due to the absence of spectral gap between
$\sigma(L)$ and the origin $\lambda=0$; see \cite{ZH}
for further discussion.)

 More precisely, it was shown in \cite{LyZ1, LyZ2}, that, similarly as in the viscous shock case \cite{ZS},
\be\label{gdrel}
D(\lambda)=\gamma \Delta(\lambda) + o(|\lambda|)
\ee
for $|\lambda|$ sufficiently small, where $\gamma$
is a constant, and $\Delta$ is a homogeneous
degree one Lopatinski determinant.

 The constant $\gamma$ measures the angle between the unstable subspace at $-\infty$ and the stable subspace at $+\infty$ for the traveling wave ODE (in this paper, equations \eqref{eq:tw1}-\eqref{eq:tw3}), that is, the \emph{transversality} of the traveling wave as a solution of the traveling wave ODE.

 The condition $\Delta(\lambda)\ne0$
for $\Re \lambda \ge 0$ and $\lambda\ne0$ is equivalent to linear stability of the corresponding inviscid shock (square wave approximation) as a solution of the hyperbolic Chapman-Jouguet equations (the Chapman-Jouguet limit is the instantaneous reaction limit, or $k = +\infty,$ of the inviscid ($b, d = 0$) limit of \eqref{model}).

Thus, similarly as in the viscous shock or relaxation case,
condition \eqref{D-condition} is equivalent to,
 \begin{equation} \label{E-condition} \left\{ \begin{aligned}
\sigma(L)\subset \{\Re \lambda\le 0\}\cup \{0\}, \\ \gamma \neq 0, \\ \mbox{$\Delta(\lambda)\ne0$
for $\Re \lambda \ge 0$ and $\lambda\ne0$}, \end{aligned} \right. \end{equation}
that is, strong spectral stability (first condition in \eqref{E-condition}), plus transversality, plus Lopatinksi stability of the associated
square-wave (Chapman--Jouguet) approximation.

 Note that \eqref{gdrel}
holds in the much more general
multidimensional case as well \cite{ZKochel, JLW}.

\br\label{gallop}
\textup{
It is shown in \cite{JLW} that, under ``standard'' assumptions
of a reaction-indepenent, ideal gas equation of state,
strong detonations are always Chapman--Jouguet stable.
Together with \eqref{gdrel}, this has the interesting consequence
that transition from viscous stability to instability
as parameters are varied must occur either by breakdown of
transversality in the traveling-wave connection, or else
by crossing of the imaginary axis of one or more nonzero complex
conjugate eigenvalue pairs, i.e., a Poincar\'e--Hopf type
bifurcation.
This agrees with physically observed
``galloping'' or ``pulsating'' instabilities;
see \cite{LyZ2, TZ1, TZ2} for further discussion.
}
\er

\begin{defi} \label{orbital} Let $X$ and $Y$ be two Banach spaces.
A traveling wave solution $\bar U$ of \emph{\eqref{RDC}} is said to be $X \to Y$
{\it nonlinearly orbitally stable} if, for any solution $\tilde U$ of \emph{\eqref{RDC}} with initial data
sufficiently close in $X$ to $\bar U,$ there exists a phase shift $\delta,$ such that
$\tilde U(\cdot,t)$ approaches $\bar U(\cdot- \delta(t))$, in $Y$ and as $t\to \infty.$
If, also, the phase $\delta(t)$ converges to a limiting value $\delta(+\infty)$,
the profile is said to be
{\it nonlinearly phase-asymptotically orbitally stable}.
\end{defi}

Using the information given by Theorem \ref{D}, we further obtain
our second main theorem, asserting that strong spectral
stability implies nonlinear stability.
This is a corollary of the more detailed,
pointwise version given in Theorem \ref{nonlin}, in which we let,
$$ \hat L^\infty := \{ f \in {\mathcal S}'(\RR), \quad (1 + |\cdot|)^{3/2} f(\cdot) \in L^\infty \}.$$
In particular, $\hat L^\infty \hookrightarrow L^1 \cap L^p,$ for all $1 \leq p \leq +\infty.$
\begin{theo}\label{nonlin} Under condition \emph{\eqref{D-condition}},
a traveling combustion wave $\bar U(x-st)$
of \emph{\eqref{model}}
is $\hat L^\infty \to L^p$ nonlinearly phase-asymptotically
orbitally stable, for $p > 1.$ More precisely, given $\bar U$ a traveling-wave solution of \emph{\eqref{model}}, given $1 \leq p \leq \infty,$ there exist $E_0 > 0,$ $C > 0,$ $\delta(\cdot)\in C^1$, and $\delta(+\infty) \in \RR,$
such that the unique solution $\tilde U$ of \emph{\eqref{RDC}} issuing from the initial datum $\bar U + U_0,$ where,
 $$ (1 + |x|)^{3/2} |U_0(x)| \leq E_0,$$
 satisfies the asymptotic estimates,
\begin{equation}
\label{LP}
\begin{aligned}
|\tilde U(x,t)-\bar U(x - \delta(t))|_{L^p}&\le C E_0
(1+t)^{-\frac{1}{2}(1-\frac{1}{p})},\\
  |\dot \delta (t)|&\le C E_0 (1+t)^{-1},\\
  |\delta(t)-\delta(+\infty)|&\le C E_0(1+t)^{-1/2}.
\end{aligned}
\end{equation}
\end{theo}

\end{subsection}


\begin{subsection}{Comments} \label{comments}

 We indicate in this section how the above theorems relate to previous mathematical results on combustion waves.

 \emph{Strong detonations} are combustion waves for which the underlying gas dynamical shock is of Lax type (see section \ref{rh-section}).

 Nonlinear stability of small-amplitude strong detonations in the small-$q$ limit was established by Li, Liu and Tan \cite{LLT}
using spectral analysis together with
Sattinger's method of weighted norms \cite{Sa} and by
Liu and Ying \cite{LYi} using energy estimates, for Majda's model.
Nonlinear stability of strong detonation in the small-$q$ limit for the related Majda-Rosales model (where $z_t$ is replaced by $z_x$ in the reaction equation), with explicit rates of convergence, was established by Li \cite{Li1}.
Nonlinear stability of small-amplitude strong detonations in the small-$q$ limit for full reactive Navier--Stokes (with Heaviside-type ignition function, and reaction-independent
equation of state) was established by Tan and Tesei in \cite{TT}
using detailed energy estimates.

 Roquejoffre and Vila \cite{RV} studied spectral stability of arbitrary amplitude strong detonations in the small-$k$ (ZND) limit, for Majda's model (in the case $d= 0 $).
Together with the weighted norm argument of \cite{LLT}, this is sufficient
to yield nonlinear stability with time-exponential rate, for exponentially
decaying initial data.

Most recently, Lyng and Zumbrun \cite{LyZ1} have shown by an elementary
perturbation argument using an abstract Evans function framework
that spectral stability of strong detonations in the small-$q$ (i.e.,
fluid-dynamical)
limit amounts to spectral stability of the underlying gas-dynamical shock.

  For strong detonations, in the Majda model case,
the spectral results of
\cite{LLT} yield \eqref{D-condition} and
thus full linearized and nonlinear stability.
 We note that, even though the weighted norm method suffices
(as pointed out in \cite{LLT}) to yield a nonlinear stability
result for strong detonations with exponentially decaying initial perturbations, our result
applies to much more general (in particular, algebraically decaying) data
and yields additional pointwise detail on solution structure.

%
%


As in the shock wave case, our approach
yields ultimately the {\it reduction} of stability analysis
to a spectral problem.
 Thus, the following result
is a consequence of Theorem \ref{nonlin},
together with the spectral analysis of Lyng and Zumbrun \cite{LyZ1}.
This greatly extends and simplifies the strongest prior result of \cite{TT},
illustrating the power of the method.

\begin{cor}\label{qcor}
Strong detonation waves of \eqref{model}
are $\hat L^\infty\to L^p$ nonlinearly orbitally stable, $p>1$,
in small-$q$ limit if and only if the limiting gas-dynamical
shock is stable; in particular, for Majda's model, they are always
stable.
\end{cor}

\br\label{sysRV}
\textup{
It would be interesting to extend the spectral analysis
of \cite{RV} to the case $d\ne 0$, which
would then imply nonlinear stability of arbitrary strength
strong detonations for Majda's model in the small-$k$ (ZND)
limit.
More interesting still would be to extend this
to the system case.
We conjecture that the proper system analog, similar to the
small-$q$ result of \cite{LyZ1}, is that stability
in the ZND limit is equivalent to gas-dynamical stability of
the component Neumann shock (see discussion, \cite{LyZ1, GS1, GS2})
together with hyperbolic (i.e., Lopatinski)
stability of the associated ZND detonation.
This should be provable by a combination of the singular perturbation
methods of \cite{PZ, FS} and ``multi-pulse'' calculations
carried out for multiple traveling-pulse solutions in models
of nerve-impulse and optical transmission.
This would recover the \cite{RV} result of stability
for the scalar, Majda model, for which the gas-dynamical shock, since
scalar, is automatically stable (see, e.g., \cite{Sa}) and the ZND
detonation may readily be calculated to be stable.
For systems, however, ZND detonations are often unstable, so that
stability in the ZND limit should not be expected.
}
\er

\medskip

 \emph{Weak detonations} are combustion waves for which the underlying gas dynamical shock is undercompressive (see section \ref{rh-section}).

 Nonlinear stability of weak detonations was obtained by Szepessy \cite{S}
for small-amplitude waves with intermediate $k,$
and by Liu and Yu \cite{LY} for arbitrary amplitude waves in the large-$k$ limit. Both these papers deal with the Majda-Rosales system.

 As far as we know, Theorem \ref{nonlin} is the
first analytical result on nonlinear stability of weak detonations for
the Majda model (more generally, the vectorial version including
reactive Navier--Stokes equations with artificial viscosity), and also for deflagrations of any type.\footnote{
In particular, we note that the weighted norm technique of \cite{LLT}
is inherently restricted to strong detonations for the scalar, Majda model;
see the discussion of \cite{TZ1}.}



 \medskip

 \emph{Deflagrations}, weak and strong, are other types of "undercompressive" combustion waves.


 It would be very interesting to determine (presumably by numerical
computations) the existence (here assumed) and
stability or instability (that is, verification of condition \eqref{D-condition}) of weak or strong deflagrations.


 Finally, note that abstract one-dimensional stability results on deflagrations
are likely not to be so physically important, since multi-d transverse
instabilities appear to play such a prominent role in their behavior \cite{B}.

As noted in \cite{M}, one-dimensional deflagration waves feature a pressure and a velocity which are locally nearly constant. Then, a ``constant density'' approximation \cite{MS}  shows that the interaction between the chemical reaction and fluid dynamics may be neglected. That is, roughly speaking, the complicated equations modeling reacting gas decouple into a part describing the fluid flow and a part describing the chemical reaction. As a result, deflagration waves are often modeled as systems of reaction--diffusion equations.

 The fact that detonations are usually approximated by reaction--convection
equations (ZND), and deflagrations by reaction--diffusion equations, reflects
 the general belief that these are dominating effects in the two different
contexts.  Our analysis here via reaction--convection--diffusion puts
both on the same footing, allowing treatment in a unified theoretical/numerical
framework, {\it investigation/validation} of these beliefs, determination
of their realms of validity.

\end{subsection}

\begin{subsection}{TO MERGE: Comments left out}



\end{subsection}


\begin{subsection}{Notes on the proof} \label{notes}

An important aspect of the Lax shock analysis is that differentiated source
leads to faster temporal decay \cite{Z3, MaZ3}.
Where $\phi\ge c_0>0$, note as in relaxation case that source in nonequibrium
mode gives faster-decaying response, at differentiated rate, and so can
be treated as in relaxation case.
Where $\phi=0$, undifferentiated source does not appear,
and so can be treated as in usual shock case.
What makes this technically feasible is that, near traveling waves,
the two regions are spatially separated, corresponding to
$x\le -M$ and $x\ge M$, $M>0$, respectively.
The intermediate regime $c_0>\phi\ge 0$ is localized
within the internal layer,  corresponding undifferentiated
source appears with exponentially decaying multiplier
$e^{-\theta|x|}$, $\theta>0$.  But, sources of the latter order appear
already in the undercompressive shock case, and can be treated by the
methods of \cite{HZ} with no change.

 So, our analysis is by interpolation between the viscous undercompressive
shock analysis of \cite{HZ, MaZ3} and the relaxation shock
analysis of \cite{MaZ1}. The new aspects of the
argument not present in the undercompressive shock case are isolated to
bound \eqref{Q}(ii), Remark \ref{emphrem},
and the new auxiliary Lemma \ref{auxconvolutions}.
In particular, no new convolution estimates were necessary, only a series of observations
having to do with the fact that undifferentiated sources appear in
a direction for which the Green function decays more rapidly, at
differentiated rate.

\begin{rems}\label{gc}
\textup{
1. As a relaxation system, \eqref{model} is somewhat degenerate,
violating the usual assumption of genuine coupling
between equilibrium and relaxation variables $u$ and $z$
associated with time-asymptotic smoothing of solutions
(see, e.g., \cite{MaZ1, Z} and references therein).
Indeed, asymptotic decoupling of the $z$-equation plays an important
role in the analysis; see Remarks \ref{blocktri} and \ref{vecform}.2.
Diffusion terms $b$, $d$ not present in standard relaxation systems
enforce smoothing directly.
}

\textup{
2. The case $d=0$ that is often considered for Majda's model requires
slightly different handling.  Absence of $z$-diffusion leads
to ``hyperbolic'' delta-function components in reactive modes
reminiscent of those encountered in \cite{MaZ1, MaZ3} in the case
of relaxation or degenerate viscosity, but with the difference that
incoming modes on side $x\ge 0$ are not time-exponentially damped.
This can be accomodated in the analysis by the introduction of an
exponentially weighted norm in the spirit of \cite{Sa}
in the $z$-component only, for $x\ge 0$, using the property of
exponential decoupling as $x\to 0$ of reactive and fluid modes.
This is a mathematical issue only;
for physical models, $d$ is strictly parabolic: $\Re \sigma(d)>0$.
}
\end{rems}

\end{subsection}


\begin{subsection}{Extension to the Navier-Stokes equations with physical viscosity}

The full reactive Navier--Stokes equations with real, or
physical viscosity may be treated
by essentially the same techniques, using the more complicated arguments
(and more detailed Green fn. bounds) developed in
\cite{MaZ3, MaZ4, Z, Raoofi, HRZ} for the treatment of viscous
shocks with real viscosity.  However, these arguments so far
are limited to the strong detonation case.  (Likewise, for technical
reasons, the viscous shock theory is so far limited for
physical viscosity to the Lax and overcompressive case.)
We leave this to a future work.

\end{subsection}


\begin{subsection}{Plan of the paper}

In Section \ref{sec:prelim}, we describe Majda's combustion model
and its vectorial generalization,
in Section \ref{profiles} the various types of traveling
wave connections that may occur,
and in Section \ref{s:evalue} the linearized eigenvalue equations
about these traveling waves.
In Sections \ref{s:evans} and \ref{resolvent},
we construct the Evans function and resolvent kernel of the Linearized
operator about the wave following the abstract framework of \cite{ZH, MaZ3},
specializing in the low-frequency regime to the special structure of
\eqref{model}
using the limiting constant-coefficient calculations
of Section \ref{s:evalue}.
In Section \ref{green}, we convert the resulting pointwise resolvent
kernel bounds to pointwise Green function bounds by stationary
phase type estimates on the Inverse Laplace transform formula,
in the process establishing Theorem \ref{D} equating linearized
and spectral (Evans) stability.
Finally, in Section \ref{nstab}, we carry out a nonlinear stability
analysis, establishing Theorem \ref{nonlin}.

\end{subsection}


\section{Preliminaries}\label{sec:prelim}

\subsection{Majda's model}\label{majda}
We begin with the scalar version of system \eqref{model}.  We assume as in \cite{LyZ2} that $f$, $\phi \in C^2$,
\begin{equation} \label{eq:flux_assumption}
f_u(u,z)>0,\quad f_{uu}(u,z)>0,
\end{equation}
and that $\phi$ is a ``bump''-type ignition function that is identically
zero for $u\leq u_i$ or $u\ge u^i$ and strictly positive for $u_i<u<u^i$.

 It is sometimes useful to rewrite \eqref{model} in the conservative form
\begin{align}
(u+qz)_t+f(u,z)_x& =u_{xx}+qdz_{xx},\label{eq:cmm1} \\
z_t & = dz_{xx}-k\phi(u)z.\label{eq:cmm2}
\end{align}

\begin{rems}\label{modeling}
\textup{
1. Note that the flux $f$, modeling equation of state,
depends on $z$, modeling the chemical constituation of the gas,
with the linear dependence loosely following the averaged equation of
state derived in \cite{CHT} for the full Euler equations.
This is important for realistic modeling of the full equations
of reacting flow;
see \cite{CHT, LyZ2} for further discussion.
For the Majda model, new qualitative phenomena emerge
for $f_z \not\le 0$ at the levels of both existence and behavior
of detonation profiles \cite{LyZ2}.
}

\textup{
2. Following \cite{M}, $\phi$ is usually taken to be
a ``step''-type function, vanishing for $u\le u_i$ and positive
for $u>u_i$.
As discussed in \cite{LyZ2}, our alternative choice of a bump-type function
is motivated by the physical parametrization of
temperature with respect to velocity $u$ in the traveling-wave
phase portrait of the ZND model.
This choice admits all the phemomena of the step-ignition case,
restricting to $u<u^i$.
In addition,  it allows for existence of weak deflagration
profiles (defined in Section \ref{profiles}),
as the step-type ignition function in general does not; see \cite{LyZ2} or
Remark \ref{defexist}.
}

\end{rems}


\subsection{Reacting Navier--Stokes equations}\label{sysprelim}

The single-species reacting Navier--Stokes equations with artificial viscosity,
written in Lagrangian coordinates, take the form
\ba \label{rNS}
\tau_t-v_x& =0,\\
v_t + p_x&= b_1 v_{xx}\\
E_t+(pv)_x& =b_2 E_{xx}+q_3k\phi(T)z, \\
z_t & = d z_{xx}-k\phi(T)z, \ea where $\tau$, $v$, $E=e+qz+v^2/2$,
$z$ denote specific volume ($\rho^{-1}$, where $\rho$ is density),
velocity, total energy density, and mass fraction of reactant,
$T=T(\tau, e)$ temperature, and $p=p(\tau, e, z)$ pressure, $k,
q_3, b_j, d>0$ constant, or
$$
q=\bp 0\\ 0\\ q_3\ep, \qquad b=\bp 0& 0 &0\\0& b_1& 0\\0& 0&
b_2\ep
$$
in \eqref{model}.
The ignition function $\phi$ is assumed to vanish identically
for $T$ below some critical ignition temperature $T_i$, and to be
strictly positive for $T$ above $T_i$.

A common choice of equation of state is the reaction-independent
gamma-law
\be\label{gammalaw}
p=\Gamma e/\tau,
\qquad
T=c^{-1} e,
\ee
where $c$ is the specific heat constant and
$\Gamma=\gamma-1$ is the Gruneisen constant.
In the thermodynamical rarified gas approximation, $\gamma>1$ is the
average over constituent particles of $\gamma=(n+2)/n$, where $n$
is the number of internal degrees of freedom of an individual particle
\cite{Ba}.

A more accurate assumption, following \cite{CHT}, is to view the
gas as a composite of unburned and burned phases with {\it different}
equations of state $T_j(\tau_j, e_j)$, $p_j(\tau_j, e_j)$, $j=1$
corresponding to the unburned and $j=2$ to the burned state, with
\ba\label{compositelaws}
\tau_1&=\tau/z,\quad \tau_2=\tau/(1-z)
\qquad \hbox{\rm (i.e., $\rho_1=z \rho$, \quad $\rho_2=(1-z)\rho$)},
\\
T_1&=T_2=T,\quad e=e_1+e_2,\quad p=p_1+p_2.\\
\ea
If both phases obey (different) gamma-law equations of state,
this leads \cite{CHT} to a gamma-law-type equation of
state \eqref{gammalaw} with reaction-dependent coefficients
\be\label{mixedeos}
c(z)= zc_1+(1-z)c_2, \qquad
\Gamma(z)= \frac{zc_1\Gamma_1+ (1-z)c_2\Gamma_2}
{zc_1+(1-z)c_2}.
\ee
This is the ``typical'' equation of state we have in mind.


\section{Traveling Waves}\label{profiles}

We consider traveling-wave solutions, i.e., solutions of the form
\[
u(x,t)=\bar u(x-st),\quad z(x,t)=\bar z(x-st),\quad s>0,
\]
of \eqref{model} that connect an unburned state $(u_+,z_+)=(u_+,1)$ to a completely burned state $(u_-,z_-)=(u_-,0)$. These are combustion waves that move from left to right leaving completely burned gas in their wake.

Thus, dropping bars for notational convenience, we find that the traveling-wave Ansatz leads, after an integration, from \eqref{model} to the system of ordinary differential equations:
\begin{align}
u'&=f(u,z)-f(u_-,z_-)-qdy-sqz-s(u-u_-), \label{eq:tw1}\\
z'&=y,\label{eq:tw2}\\
y'&=d^{-1}\big(-sy+k\phi(u)z\big), \label{eq:tw3}
\end{align}
where we have used $y:=z'$ to write the system in first order and $'$ denotes differentiation with respect to the variable $\xi:=x-st$.
We assume that the end states are such that
\begin{equation} \label{eq:u_endstates}
u_- \in [u_i,u^i],
\quad  u_+\not \in [u_i,u^i],
\end{equation}
so that
\begin{equation}\label{eq:phi_endstates}
\phi(u_-)>0,\quad \phi(u_+)=0,\quad \phi'(u_+)=0.
\end{equation}
Equation~\eqref{eq:u_endstates} has the physical interpretation that the unburned end state is below ignition temperature
so that the there is no chemical reaction on the unburnt side.

\subsection{Rankine--Hugoniot conditions} \label{rh-section}
A necessary condition for the existence of a connection is that the end states at $\pm\infty$ be rest points of the traveling-wave equation.
This leads to the Rankine-Hugoniot condition
\begin{equation}
f(u_+,z_+)-f(u_-,z_-)=sq+s(u_+-u_-),\label{eq:rh}
\tag{RH}
\end{equation}
together with the requirements that
\begin{equation}
y_\pm=z'_\pm=0
\end{equation}
and (justifying assumptions \eqref{eq:phi_endstates})
\begin{equation}
\phi(u_\pm)z_\pm=0.
\end{equation}

Restricting now to Majda's model, write
$\ha_\pm:=f_u(u_\pm,z_\pm)$ and
$\beta_\pm:=f_z(u_\pm,z_\pm)$.
 Then, the traveling-wave profile is
said to be a \emph{strong detonation} if
\begin{equation}
\ha_->s>\ha_+.
\label{eq:strongdet}
\end{equation}
It is said to be a \emph{weak detonation} if
\begin{equation}
s> \ha_-, \, \ha_+.
\label{eq:weakdet}
\end{equation}
Similarly, it is said to be a {\it weak deflagration} if
\begin{equation}
\label{eq:weakdefl}
\ha_-, \ha_+>s.
\end{equation}
It is said to be a {\it strong deflagration} if
\begin{equation}
\ha_+>s> \ha_-.
\label{eq:strongdefl}
\end{equation}
Degenerate
profiles for which the inequalities are nonstrict
are called {\it Chapman--Jouguet} detonations or deflagrations
and lie on the boundary between weak and strong branches.

\br\label{defexist}
\textup{
For $f$ independent of $z$, we find by \eqref{eq:flux_assumption}
that detonations correspond to case $u_->u_i\ge u_+$,
deflagrations to case $u_-\le u^i <u_+$.  In particular,
deflagrations cannot occur for a step-type ignition function,
for which $u^i=+\infty$.
}
\er

A routine modification of Propositions 2.1 and 2.2, \cite{LyZ2}, accounting
for $z$-dependence of $f$, yields the following description of solutions
of \eqref{eq:rh}.

\begin{prop}\label{RHexist}
For fixed $u_+$, suppose that $f(u_+, 1)<f(u_++q,0)$.
Then, there exist $s^*(u_+)<s_*(u_+)$ such that
(i) for $s>s_*$ there
exist two states $u_->u_+$ for which \eqref{eq:rh}
(but not necessarily \eqref{eq:phi_endstates})
is satisfied (weak and strong detonation),
for $s=s_*$ there exists one (Chapman--Jouguet detonation),
and for $s< s_*$, there exist none
no solutions $u_->u_+$. (ii)
For $s<s^*$, there exist two states $u_-<u_+$ for which \eqref{eq:rh}
is satisfied (weak and strong deflagration),
for $s=s^*$, there exists one (Chapman--Jouguet deflagration),
and for $s> s^*$, there exist none.
If $f(u_+, 1) \ge f(u_++q,0)$, on the other hand,
then for each $s$ there exists
at most one solution $u_->u^+$ and one solution $u_-<u^+$,
(strong detonation and strong deflagration, respectively).\footnote{
Typically, both; in particular, if $f$ grows superlinearly in $|u|$, then
both solutions exist for each $s$.}
\end{prop}

\begin{rem}\label{cjrem}
\textup{
The case $f(u_+, 1)<f(u_++q,0)$ is essentially identical to that
of the reaction-independent case discussed in \cite{LyZ2},
for which $f(u_+)<f(u_++q)$ by \eqref{eq:flux_assumption}.
More generally, $f_z\le 0$ is sufficient for $f(u_+, 1)<f(u_++q,0)$.
For further discussion of the \eqref{eq:rh}
problem,
see \cite{LyZ2}: in particular the Chapman--Jouguet diagrams of Figure 1.
} \end{rem}

\subsection{The connection problem}
Linearizing \eqref{eq:tw1}--\eqref{eq:tw3} around the state $(u_-,z_-,y_-)$, we find the constant-coefficient system of ordinary differential equations
\begin{equation}
\begin{pmatrix} u \\ z \\ y \end{pmatrix}'=
\begin{pmatrix}
\ha_--s & b_--sq & -qd \\
0 & 0 & 1 \\
0 & kd^{-1}\phi(u_-) & -sd^{-1}
\end{pmatrix}\begin{pmatrix} u \\ z \\ y \end{pmatrix}.
\label{eq:minus_lin_twode}
\end{equation}
For strong detonations,
the coefficient matrix in \eqref{eq:minus_lin_twode} is easily seen to have two positive eigenvalues and one negative eigenvalue. Thus, there is a two-dimensional unstable manifold at $(u_-,0, y_-)$. Similarly, we note that the linearization of \eqref{eq:tw1}--\eqref{eq:tw3} about $(u_+,z_+,y_+)$ is
\begin{equation}
\begin{pmatrix} u \\ z \\ y \end{pmatrix}'=
\begin{pmatrix}
\ha_+-s & b_+-sq & -qd \\
0 & 0 & 1 \\
0 & 0 & -sd^{-1}
\end{pmatrix}\begin{pmatrix} u \\ z \\ y \end{pmatrix}.
\label{eq:plus_lin_twode}
\end{equation}
By the block-triangular structure, it is easy to see that there are two negative eigenvalues  and one zero eigenvalue. It is easy to see that the center manifold is a line of equilibria, so plays no orbit may approach the rest point $(u_+,z_+,y_+)$ along the center manifold. Thus, for connections, the important structure is the two-dimensional stable manifold at $(u_+,1, y_+)$. Counting dimensions, we see that a connection corresponds to the intersection of two two-dimensional manifolds in $\mathbb{R}^3$.
In particular, it generically persists as a unique, transverse intersection,
under variations in parameters such as $u_\pm$, $s$ consistent with
\eqref{eq:rh}.
See \cite{LyZ2} for a discussion of this situation in the case $d=0$.

Similarly, for weak detonations, there is
a one-dimensional unstable manifold at $(u_-,0, y_-)$
and a two-dimensional stable manifold at $(u_+,0, y_+)$.
Thus, connections are typically {\it codimension one} in the
set of \eqref{eq:rh} compatible parameters, in contrast to the strong
detonation case.
See \cite{LyZ2} in the case $d=0$.
This situation is analogous to that of a Lax-type shock in the nonreactive case.
For weak deflagrations, there is a
a two-dimensional unstable manifold at $(u_-,0, y_-)$
and a one-dimensional stable manifold at $(u_+,0, y_+)$.
Thus, connections are again generically {\it codimension one} in the
set of \eqref{eq:rh} compatible parameters.
See \cite{LyZ2} in the case $d=0$.
For strong deflagrations, there is a
a one-dimensional unstable manifold at $(u_-,0, y_-)$
and a one-dimensional stable manifold at $(u_+,0, y_+)$,
and connections are {\it codimension two}.

In every case, we have by the discussion above:

\begin{lem}\label{l:expdecay}
Traveling-wave profiles $(\bar u, \bar z)$ corresponding to weak or strong
detonations or deflagrations satisfy
\be\label{e:expdecay}
|(d/dx)^k \Big((\bar u, \bar z)- (u,z)_\pm\Big)|\le Ce^{-\theta |x|},
\qquad
x\gtrless 0, \quad 0\le k\le 3.
\ee
\end{lem}

\begin{proof}
Standard ODE estimates for stable and unstable manifolds.
\end{proof}

\br\label{class}
\textup{
Weak detonations and deflagrations are analogous to undercompressive
shocks in the nonreactive case, with strong deflagrations undercompressive
of degree two; see \cite{Z} for a discussion of shock classification.
In the case $d=0$, it can
be demonstrated that weak detonation connections do occur in some
cases, but deflagration connections (weak or strong) of the type we have
described do not \cite{LyZ2}.\footnote{More precisely, weak deflagration
profiles exist
only in the degenerate case $u_+=u^i$, with $u(x)$ converging
to $u_-$ as $x\to +\infty$ at subalgebraic rate;
strong deflagration profiles do not exist in any case \cite{LyZ2}.}
It is an interesting question whether or not they occur for $d\ne 0$,
or, more generally, for the full, reactive Navier--Stokes equations
\cite{LyZ1, LyZ2}.
}
\er

\subsection{The system case}\label{sysprofile}
Under the further assumption of (asymptotic) {\it dissipativity},
\be\label{diss}
\Re \sigma( \ha i \xi -\xi^2 b)_\pm\le \frac{-\theta |\xi|^2}
{1+|\xi|^2},
\qquad \ha_\pm:=(\partial f/\partial u)(U_\pm),
\ee
$\theta>0$, for all $\xi\in \RR$ (standard for systems \cite{Z}),
it is readily verified using the above-mentioned block-triangular
decomposition of limiting systems
into fluid and reactive blocks together with fluid
dynamical results of Majda and Pego \cite{MP}, that the main results
of this section carry over to the system case, substituting
for definitions \eqref{eq:strongdet}--\eqref{eq:strongdefl} the system
versions (for right-going waves, $s>0$)
\begin{equation}
\ha_n^->s>\ha_n^+,
\quad
s>\ha_j^\pm, \, j\ne n,
\qquad
\hbox{\rm (strong detonation)}
\label{eq:sysstrongdet}
\end{equation}
\begin{equation}
s> \ha_n^-, \, \ha_n^+,
\quad
s>\ha_j^\pm, \, j\ne n,
\qquad
\hbox{\rm (weak detonation)}
\label{eq:sysweakdet}
\end{equation}
\begin{equation}
\label{eq:sysweakdefl}
\ha_n^-, \ha_n^+>s,
\quad
s>\ha_j^\pm, \, j\ne n,
\qquad
\hbox{\rm (weak deflagration)}
\end{equation}
and
\begin{equation}
\ha_n^+>s> \ha_n^-.
\quad
s>\ha_j^\pm, \, j\ne n,
\qquad
\hbox{\rm (strong deflagration)}
\label{eq:sysstrongdefl}
\end{equation}
where $\ha_1^\pm< \cdots < \ha_n^\pm$ denote the eigenvalues of
$\ha_\pm:=f_u(u_\pm,z_\pm)$.

In particular, connections if they exist satisfy \eqref{eq:rh}, and
if they are transverse are generically of codimension zero, one, one, two,
respectively, in the set of \eqref{eq:rh}-compatible parameters.
Further, the connecting profile satisfies \eqref{e:expdecay}, converging
exponentially to its endstates $U_\pm$ as $x\to \pm \infty$.
See \cite{ZKochel}, Appendix A for further details.
Condition \eqref{diss} holds trivially for identity viscosity $b= I$,
and holds also for the physical (semidefinite) viscosity
of the reacting Navier--Stokes equations \cite{MaZ4, Z},
the two main cases we have in mind.
For simplicity of notation, we assume also that $\sigma(d)$ is
semisimple, so that $d$ has a full set of eigenvectors.

Likewise, there is a simple analogy to the Chapman--Jouget
analysis of Section \ref{rh-section}, and, for the typical
mixed gamma-law equation of state \eqref{mixedeos}
of Proposition \ref{RHexist}.
For, rearranging \eqref{eq:rh} in the case of the
reacting Navier--Stokes equations \eqref{rNS},
we obtain \cite{LyZ1} from the first equation
that $(v_+-v_-)=-s(\tau_+-\tau_-)$,
from the second that $(\tau, p)_\pm$ lie on the Rayleigh line
\be\label{Ray}
p_+- p_-=-s^2 (\tau_+-\tau_-),
\ee
and from the third the shifted Hugoniot curve
\be\label{Hug}
(e_+-e_-) + q=(-1/2)(p_++p_-)(\tau_+-\tau_-).
\ee

Thus, fixing $(\tau, v, E)_+$,
viewing \eqref{Hug} as determining a ``burned'' pressure law
\be\label{pressurelaw}
p_-=P_-(\tau_-),
\ee
and assuming that $e_+$ can be recovered from $\tau_-$, $p_-$
through inversion of $p_-=p(\tau_-, e_-, 0)$,
we find that the allowable states
$(\tau, v, E)_-$ are determined as the intersection in the $\tau-p$
plane of line \eqref{Ray} with curve \eqref{pressurelaw}, similarly as in
the scalar (Majda's model) case.
Moreover, there are two distinct solution structures, according
as
\be\label{solcrit}
p_+< P_-(\tau_+)
\ee
(standard: for $P_-$ convex, pairs of weak/strong detonations, deflagrations
as in the scalar case) or
the reverse (nonstandard: for $P_-$ convex,
single strong detonation, deflagration).

For the typical equation of state \eqref{gammalaw}--\eqref{mixedeos},
it is readily calculated that
\be\label{mixedpress}
P_-(\tau)=
\frac{ \big( \tau_+/\Gamma_1 - (1/2)(\tau - \tau_+)\big)p_+ + q}
{ \tau/\Gamma_2 - (1/2)(\tau - \tau_+)},
\ee
whence \eqref{solcrit} reduces to
\be\label{idealcrit}
p_+(1-\Gamma_2/\Gamma_1)< q\Gamma_2/\tau_+.
\ee
From \eqref{idealcrit}, we see that
$\Gamma_1\le \Gamma_2$ (in particular, including the reaction-independent
case) implies a standard \eqref{eq:rh} solution structure.
Roughly speaking, this corresponds to a reaction in which
complicated compounds break up into simpler components,
so that $n$ decreases and $\Gamma=2/n$ increases,
the reverse situation to a reaction in which simple
components combine into more complicated molecules.
We conjecture, by analogy with the scalar case, that
for \eqref{mixedeos}, condition \eqref{idealcrit}
equivalent to $\Gamma_z \le 0$ implies further a standard connection
structure, at least in the ZND limit (for which there is a close
connection to Majda's model \cite{LyZ1, LyZ2, GS1}).
However, even for $\Gamma_1>\Gamma_2$, the standard \eqref{eq:rh}
solution structure is recovered for $q$ sufficiently large.

\begin{rems}\label{sysconnection}
\textup{
1. Dissipativity, \eqref{diss}, implies in particular
hyperbolicity of the first-order convection terms, i.e.,
that $\alpha$ has real, semisimple eigenvalues.
We use this freely below.
}

\textup{
2. Existence of detonation connections for the full, reacting
Navier--Stokes equations
has been studied by Gasser and Szmolyan \cite{GS1, GS2}
using geometric singular perturbation techniques in the ZND
($b$, $d\to 0$) limit, and by Gardner \cite{G} using Conley
index methods.
}
\end{rems}

\section{The eigenvalue equation}\label{s:evalue}

Suppose $\big( \bar u(x-st),\bar z(x-st)\big)$ is a
traveling-wave profile of \eqref{model} as described above. We now begin to investigate the stability of such an object.
The linearized equations about $(\bar u, \bar z)$ in moving coordinates
$\tilde x =  x-st$, are, dropping tildes,
\begin{align}
&u_t-q(k\phi'(\bar u)u\bar z+k\phi(\bar u)z)+(\alpha u)_x+(\beta z)_x=u_{xx},\label{eq:mmlin1}\\
&z_t-sz_x=-k\phi'(\bar u)u\bar z-k\phi(\bar u)z+dz_{xx},\label{eq:mmlin2}
\end{align}
where $\alpha:=f_u(\bar u,\bar z)-s$, $\beta:=f_z(\bar u,\bar z)$, and $u$, $z$ now denote perturbations. The eigenvalue equations corresponding to this linear system are thus
\begin{align}
&u''=\big(\lambda u-q(k\phi'(\bar u)u\bar z+k\phi(\bar u)z)+(\alpha u)'+(\beta z)'\big),\label{eq:eval1} \\
&z''=d^{-1}\big(\lambda z-sz'+k\phi'(\bar u)u\bar z+k\phi(\bar u)z\big),\label{eq:eval2}
\end{align}
Alternatively, upon substituting $dz''-\lambda z+sz'=k\phi'(\bar u)u\bar z + k\phi(\bar u)z$ from \eqref{eq:eval2} into \eqref{eq:eval1}, we can rewrite \eqref{eq:eval1} as
\begin{equation}
u''=\big(\lambda(u+qz)-sqz'-qdz''+(\alpha u)'+(\beta z)'\big).\label{eq:alteval1}
\end{equation}
Compare this with the remark at the end of Section~\ref{sec:prelim}.
We write \eqref{eq:eval1}--\eqref{eq:eval2} as a first-order system. To do so, we define
$W:=(u,z,u',z')^\trans$, so that \eqref{eq:eval1}--\eqref{eq:eval2} becomes
\be\label{eval}
W'=\mA(x,\lambda)W
\ee
where the coefficient matrix is
\be
\label{eq:Amatrix}
\mA(x,\lambda) =
\bp
0 & 0 & 1 & 0  \\
0 & 0 & 0 & 1 \\
\lambda + \alpha'-qk\phi'(\bar u)\bar z & \beta' -qk\phi(\bar u) & \alpha & \beta \\
d^{-1}k\phi'(\bar u) \bar z & d^{-1}\lambda + d^{-1}k\phi(\bar u) & 0 & -sd^{-1}
\ep.
\ee
System \eqref{eval} has a limiting constant-coefficient structure, i.e., the coefficient matrix
has limits as $x\to\pm\infty$. That is,
\[
\mA(x,\lambda)\to\mA_\pm(\lambda)\quad\text{as}\quad x\to\pm\infty,
\]
and the limiting matrices are given by
\be
\label{eq:AAlimits}
\mA_\pm(\lambda)=
\bp
0 & 0 & 1 & 0 \\
0 & 0 & 0 & 1 \\
\lambda & -kq\phi(u_\pm) & \alpha_\pm & \beta_\pm \\
0 & d^{-1}(\lambda+k\phi(u_\pm)) & 0 & sd^{-1}
\ep.
\ee

\br\label{blocktri}
\textup{Here, we coordinatize as $W=(u,z,u',z')^\trans$ following the general
scheme of \cite{ZH, MaZ3} as this is what we shall need below to establish the pointwise bounds. However, it is sometimes helpful to see the fluid/reaction structure in the system. To see this at the level of the eigenvalue ODEs, we write
\be
\hat W:=(u,u',z,z')^\trans,
\label{eq:altw}
\ee
separating the fluid $(u,u')$ and reaction $(z,z')$ quantities. In this labeling scheme, the eigenvalue ODE~\eqref{eval} becomes
\be\label{eq:alteval}
\hat W'=\hat\mA(x,\lambda)\hat W
\ee
where the matrix $\hat \mA$ can easily be obtained from \eqref{eq:Amatrix} by appropriately swapping entries, and
\be
\label{eq:altlimits}
\hat\mA_\pm(\lambda)=\bp
0 & 1 & 0 & 0 \\
\lambda & \alpha_\pm & -kq\phi(u_\pm) & \beta_\pm \\
0 & 0 & 0 & 1 \\
0 & 0 & d^{-1}(\lambda+k\phi(u_\pm)) & -sd^{-1}
\ep.
\ee
In particular, the upper block-triangular structure of $\hat\mA_\pm$ will be useful in the constant-coefficient analysis below.}
\er

\subsection{Constant--coefficient analysis}\label{cc}

We now examine the constant-coefficient limiting systems $W'=\mA_\pm(\lambda)W$ (or equivalently $\hat W'=\hat\mA_\pm(\lambda)\hat W$) in some detail. From the upper block-triangular structure in \eqref{eq:altlimits}, it is quite straightforward to compute eigenvalues; they are simply the eigenvalues of the diagonal blocks. From the upper left-hand ``fluid'' block, we obtain the fluid eigenvalues
\be
\label{eq:fluideval}
\mu^\pm_f=\frac{\alpha_\pm\pm\sqrt{\alpha_\pm^2+4\lambda}}{2},
\ee
while the lower right-hand ``reaction'' block contributes eigenvalues of form
\be
\label{eq:reactevalplus}
\mu^+_r=\frac{-sd^{-1}\pm\sqrt{s^2d^{-2}+4d^{-1}\lambda}}{2}\quad\text{from}\quad\hat\mA_+,
\ee
and
\be
\label{eq:reactevalminus}
\mu^-_r=\frac{-sd^{-1}\pm\sqrt{s^2d^{-2}+4(d^{-1}\lambda+d^{-1}k\phi(u_-))}}{2}\quad\text{from}\quad\hat\mA_-.
\ee
The corresponding eigenvectors also have structure inherited from the block-triangular nature of the limiting matrices $\hat\mA_\pm$. In particular, as long as the fluid and reaction eigenvalues remain distinct --- as our calculations below show they are for small $\lambda$, the corresponding eigenvectors take the form
\be\label{eq:cc_evecs}
\hat v_f^\pm=\bp 1 \\ \mu_f^\pm \\ 0 \\ 0 \ep
,\quad\text{and}\quad
\hat v_r^\pm=\bp M_\pm^{-1} \bp 1 \\ \mu_r^\pm\ep \\ 1 \\ \mu_r^\pm\ep,
\ee
where
\be
M_+(\lambda):=\left[\bp 0 & 1 \\ \lambda & \alpha_+\ep -\mu_r^+I\right]^{-1}
\bp
0 & 0 \\
0 & \beta_+
\ep
\ee
and
\be
M_-(\lambda):=\left[\bp 0 & 1 \\ \lambda & \alpha_-\ep -\mu_r^-I\right]^{-1}
\bp
0 & 0 \\
-qk\phi(u_-) & \beta_-
\ep.
\ee

We also record here Taylor series expansions of those eigenvalues
of the limiting systems which become zero at $\lambda=0$, the
so-called slow modes. These are \be\label{eq:mu_r_taylor}
\mu^+_r(\lambda)=\frac{1}{s}\lambda
-\frac{2d}{s^3}\lambda^2+\cdots \ee and
\be\label{eq:mu_f_taylor}
\mu_f^\pm(\lambda)=\pm\frac{1}{\alpha_\pm}\lambda\mp\frac{1}{\alpha_\pm^3}\lambda^2+\cdots,
\ee
with associated (right) eigenvectors, now written in $(u,z,u',z')$
coordinates, \be\label{eq:v_r_taylor} v^+_r(\lambda)=\bp R_r^+ \\
0 \ep + \cdots \ee and \be\label{eq:v_f_taylor} v_f^\pm(\lambda)=
\bp R_f^\pm \\ 0 \ep + \cdots, \ee where limiting fluid modes
\be\label{fluidform} R_f^\pm= \bp *\\0\ep \ee have vanishing
$z$-component.

We shall also have need of the adjoint eigenvalue equation
\ba\label{adjeval}
\tilde u''& = \lambda^*\tilde u -qk\phi'(\bar u)\bar z \tilde u + k \phi'(\bar u)\bar z\tilde z-\alpha \tilde u', \\
\tilde z''&=d^{-1}\big(\lambda^*-qk\phi(\bar u)\tilde u+k\phi(\bar
u) \tilde z-\beta\tilde u'+s\tilde z'\big) \ea associated with
\eqref{eq:eval1}--\eqref{eq:eval2}, where $\lambda^*$ denotes
complex conjugate. Writing as a first-order system $\tilde
W'=\tilde \mA(x, \lambda^*)\tilde W$, $\tilde W=(\tilde u, \tilde
z, \tilde u', \tilde z')^\trans$, and studying the eigenvalues
$\tilde \mu_j^\pm$ and eigenvectors $\tilde v_j^\pm$ of the
limiting, coefficient-matrix $\tilde \mA_\pm(\lambda^*)$, we have
by duality that $\tilde \mu_j^\pm=(\mu_j^\pm)^*$, while a brief
calculation yields Taylor expansion \be\label{dualode} \tilde
v_f^\pm(\lambda)= \bp L_f^\pm \\ 0 \ep + \cdots, \qquad \tilde
v_r^\pm(\lambda)= \bp L_r^\pm \\ 0 \ep + \cdots \ee at
$\lambda=0$, where limiting left fluid modes
 $L_f^-= c(1
,q)^\trans$
on the minus infinity (reactive) side satisfy
\be\label{dfluidform} L_f^- \perp (-q, 1)^\trans. \ee
(See Remark \ref{vecform}.1 below for further discussion, and
extension to the system case.)

\begin{defi}\label{consplit}
\textup{
The \emph{domain of consistent splitting} for an ODE $W'=\mA(x,\lambda)W$
with asymptotically constant coefficients
is the set of $\lambda\in\CC$ such that
\begin{enumerate}
\item[(i)] the limiting matrices $\mA_\pm(\lambda)$ are hyperbolic, i.e., they have no center subspace, and
\item[(ii)] the dimensions of the stable (unstable) subspaces $S^+(\lambda)$ and $S^-(\lambda)$ ($U^+(\lambda)$ and $U^-(\lambda)$) are the same.
\end{enumerate}
}
\end{defi}

\begin{lem}\label{lem:rhp}
The set $\{\lambda\in\CC\; :\; \Re\lambda>0\}$ is a subset of the domain of consistent splitting.
\end{lem}
\begin{proof}
Immediate by inspection of formulas \eqref{eq:fluideval}--\eqref{eq:reactevalminus}.
\end{proof}
In fact, more can be said. The limiting matrices $\hat\mA_\pm$ fail to be hyperbolic if
\[
\det(\hat\mA_\pm(\lambda)-i\xi)=0
\]
for some $\xi\in\RR$. But, again using the upper block-triangular structure of $\hat\mA_\pm$, it is clear that this determinant vanishes if and only if the determinant of one of the diagonal blocks vanishes. This leads to the following four \emph{dispersion curves} in the complex $\lambda$-plane,
\begin{align*}
\lambda_f^\pm(\xi)&=-\xi^2-i\xi\alpha_\pm, \\
\lambda_r^+(\xi)& = d(-\xi^2+isd^{-1}),\\
\lambda_r^-(\xi)&= d(-\xi^2+isd^{-1})-\phi(u_-).
\end{align*}
These curves, parabolae opening into the left complex half plane, can be used to describe the boundary of the domain of consistent splitting. In particular, we note that Lemma~\ref{lem:rhp} gives that the open right half complex plane is contained in the domain of consistent splitting, and, varying $\lambda$ from right to left, we find that $\lambda$ cannot leave the domain unless we cross one of these four dispersion curves.
Thus, the component of the domain of consistent splitting which contains $+\infty$ contains a set of form
\[
\Omega_\eta:=
\{\lambda: \Re \lambda >
\max\{-\eta_1 |\Im \lambda|, \eta_2 |\Im \lambda|^2\} \},
\qquad
\eta_j>0.
\]
See \eqref{containment} and Lemma~\ref{frozen} below.

\subsection{The system case}\label{syseval}
The above calculations extend in straightforward fashion to the system
case, substituting block matrix computations for the scalar computations
above.
In particular, the block-triangular structure of
\eqref{eq:altlimits} is maintained, reducing the computation of
constant-coefficient modes to a computation on the upper lefthand
diagonal fluid block that is exactly the viscous shock computation done in
\cite{ZH, MaZ3}, and a computation on the lower righthand
diagonal reaction block that on the plus infinity side is a particularly
simple (scalar convection) version of the same viscous shock computation
and on the minus infinity side consists of fast modes that need not be
resolved.

Specifically, the Taylor expansions of slow modes at $\lambda=0$
become \be\label{syseq:mu_r_taylor}
\mu^+_{r,i}(\lambda)=\frac{1}{s}\lambda
-\frac{2d_j^+}{s^3}\lambda^2+\cdots, \ee $i=1, \dots, m$, $z\in
\RR^m$,
\be\label{syseq:mu_f_taylor}
\mu_{f,j}^\pm(\lambda)=\pm\frac{1}{\alpha_j^\pm}\lambda\mp\frac{b_j^\pm}{(\alpha_j^\pm)^3}\lambda^2+\cdots,
\ee
$j=1, \dots, n$, $u\in \RR^n$, \be\label{syseq:v_r_taylor}
v^+_{r,i}(\lambda)=\bp R_{r,i}^+ \\ 0 \ep + \cdots \qquad
R_{r,i}^\pm= \bp *\\r_{r,i}\ep, \ee and
\be\label{syseq:v_f_taylor} v_f^\pm(\lambda)= \bp R_{f,j}^\pm \\ 0
\ep + \cdots, \qquad R_{f,j}^\pm= \bp r_{f,j}\\0\ep, \ee
\be\label{syseq:tildev_f_taylor} \tilde v_f^\pm(\lambda)= \bp
L_{f,j}^\pm \\0 \ep + \cdots, \qquad L_{f,j}^+= c \bp l_{f,j}^+\\
0\ep, \qquad L_{f,j}^-= c \bp l_{f,j}^-\\ q^\trans l_{f,j}^-\ep,
\ee where $l_{r,i}^+$ and $r_{r,i}$ are left and right
eigenvectors of $d$ (now matrix-valued) and $\alpha_{f,j}^\pm$,
$l_{f,j}^\pm$, $r_{f,j}^\pm$ are the eigenvalues and left and
right eigenvectors of $\partial f/\partial u(U_\pm)$,
$b_j^\pm=(l_{f,j}br_{f,j})_\pm$ ({\it strictly positive}, by
dissipativity assumption \eqref{diss}), and
$d_i^+=(l_{r,i}dr_{r,i})_+$. Note that we again have vanishing of
the $z$-component of asymptotic fluid modes $R_{f,j}^\pm$, as well
as the key orthogonality relation \be\label{sysorthog}
L_{f,j}^-\perp \bp -q \\  I_r \ep. \ee

\begin{rems}\label{vecform}
\textup{
1. The structural relations \eqref{fluidform}, \eqref{dfluidform} and their
system analogs \eqref{syseq:v_f_taylor}--\eqref{sysorthog} for
the asymptotic modes,
play an important role in our analysis;
see Remarks \ref{LFrems}.1, \ref{zdecay} and the proofs of
Proposition \ref{greenbounds}, Lemma \ref{auxconvolutions},
and Proposition \ref{pwnonlin}.
Indeed, this is essentially the only structure that we use, other
than the existence of Taylor expansions of slow modes at $\lambda=0$.
These may be verified easily by substituting into the limiting,
constant-coefficient eigenvalue systems
\ba\label{syscceval}
u''&=\big(\lambda u-qk\phi_\pm z+\alpha_\pm u'+\beta_\pm z'\big),\\
&z''=d^{-1}\big(\lambda z-sz'+k\phi_\pm z\big)
\ea
and their dual, adjoint systems
the Ansatze $U=e^{\mu x}R$, $\tilde U=e^{\mu^* x}L$, respectively, to
obtain
characteristic equations
$$
\bp
\mu^2 -\alpha_\pm \mu -\lambda I_n & qk\phi_\pm -\beta_\pm \mu\\
0 & d\mu^2 - \lambda I_r +s \mu - k\phi_\pm) \\
\ep R=0
$$
and
$$
L^*
\bp
\mu^2 +\alpha_\pm \mu -\lambda I_n & qk\phi_\pm +\beta_\pm \mu\\
0 & d\mu^2 -\lambda I_r +s \mu - k\phi_\pm) \\
\ep
=0,
$$
respectively, which, setting $\lambda=\mu=0$ to obtain slow mode behavior
at $\lambda=0$, reduce to
$$
\bp
0 & qk\phi_\pm \\
0 & - k\phi_\pm \\
\ep
R
=0,
\qquad
L^*
\bp
0 & qk\phi_\pm \\
0 & - k\phi_\pm \\
\ep
=0,
$$
from which \eqref{syseq:v_f_taylor}--\eqref{sysorthog} are evident.
}

\textup{
2. The containment $\Omega_\eta\subset \Lambda$ likewise carries
over to the system case, by dissipativity assumption
\eqref{diss} and Lemma \ref{frozen} below.
}
\end{rems}

\section{The Evans function}\label{s:evans}

We now construct the Evans function following the abstract
framework of \cite{MaZ3, Z}.
For historical origins of the Evans function, see \cite{AGJ, PW, GZ}
and references therein.

\subsection{The Conjugation Lemma}\label{s:conjugation}
We first recall a central result connecting variable- and
constant-coefficient ODE.
Consider a general family of first-order ODE
\begin{equation}
\label{gfirstorder}
W'-{\mathbb A}(x, \lambda)W=0
\end{equation}
of the form \eqref{eval},
indexed by a spectral parameter $\lambda \in \Omega \subset \CC$, where
$W\in \CC^N$, $x\in \RR$ and ``$'$'' denotes $d/dx$,
assuming (cf. Lemma \ref{l:expdecay})
\medbreak
\textup{(h0) }
Coefficient ${\mathbb A}(\cdot,\lambda)$, considered
as a function from $\Omega$ into
$C^0(x)$ is analytic in $\lambda$.
Moreover, ${\mathbb A}(\cdot, \lambda)$ approaches
exponentially to limits $\mA_\pm$ as $x\to \pm \infty$,
with uniform exponential decay estimates
\begin{equation}
\label{expdecay2}
|(\partial/\partial x)^k(\mA- \mA_\pm)|
\le C_1e^{-\theta|x|/C_2}, \,
\quad
\text{\rm for } x\gtrless 0, \, 0\le k\le K,
\end{equation}
$C_j$, $\theta>0$,
on compact subsets of $\Omega $.
\medbreak

The following asymptotic ODE result generalizes
the Gap Lemma of \cite{GZ}; for a proof, see, e.g., \cite{MZ1, MaZ3, Z}.

\begin{prop}[{The Conjugation Lemma \cite{MZ1}}]
\label{conjugation}
Given (h0), there exist locally to any given $\lambda_0\in \Omega $
invertible linear transformations $P_+(x,\lambda)=I+\Theta_+(x,\lambda)$ and
$P_-(x,\lambda) =I+\Theta_-(x,\lambda)$ defined
on $x\ge 0$ and $x\le 0$, respectively,
$\Phi_\pm$ analytic in $\lambda$ as functions from $\Omega$
to $C^0 [0,\pm\infty)$, such that:
\medbreak
(i)
For any fixed $0<\btheta<\theta$ and $0\le k\le K+1$,  $j\ge 0$,
\begin{equation}
\label{Pdecay}
|(\partial/\partial \lambda)^j(\partial/\partial x)^k
\Theta_\pm |\le C(j) C_1 C_2 e^{-\theta |x|/C_2}
\quad
\text{\rm for } x\gtrless 0.
\end{equation}
\smallbreak
(ii)  The change of coordinates $W=:P_\pm Z$
reduces \eqref{gfirstorder} to
\begin{equation} \label{glimit}
Z'-\mA_\pm Z = 0
\quad
\text{\rm for } x\gtrless 0.
\end{equation}
\end{prop}

\subsection{Normal modes}\label{normal}
Using Proposition \ref{conjugation},
we next construct normal modes for \eqref{gfirstorder}.
Recall the domain of consistent splitting defined in Definition \ref{consplit},
Section \ref{cc}.

\begin{lem}\label{Kat}
On any simply connected subset of the domain of consistent
splitting $\Lambda$, there exist analytic bases
$\{V_{1}, \dots, V_{k}\}^\pm$ and $\{V_{k+1}, \dots, V_{N}\}^\pm$
for the subspaces $S_\pm$ and $U_\pm$ defined in Definition \ref{consplit}.
\end{lem}

\begin{proof}
By spectral separation of $U_\pm$, $S_\pm$, the associated
(group) eigenprojections are analytic.
The existence of analytic bases
then follows by a standard result of Kato; see \cite{Kat}, pp. 99--102.
\end{proof}

By Lemma \ref{conjugation}, on the domain of consistent
splitting, the subspaces
\be\label{calS}
\CalS^+=\Span\{W_1^+,\dots,   W_k^+\}
:=\Span\{P_+V_1^+,\dots,  P_+ V_k^+\}
\ee
and
\be\label{calU}
\CalU^-:=\Span\{  W_{k+1}^-, \dots,  W_{N}^-\}
:=\Span\{ P_- V_{k+1}^-, \dots, P_- V_{N}^-\}
\ee
uniquely determine the stable manifold as $x\to +\infty$
and the unstable manifold as $x\to -\infty$ of \eqref{gfirstorder},
defined as the manifolds of solutions decaying as $x\to \pm\infty$,
respectively, independent of the choice of $P_\pm$.
More generally, $W_j^\pm:=(P V_j)_\pm$, $j=1, \dots, N$
are called  {\it normal modes} for \eqref{gfirstorder}.

In the context of Majda's model (more generally, the system
analog HERE),
$V_j^\pm$ are comprised of the vectors
described in Section \ref{cc} (resp. \ref{syseval}),
of which the {\it slow modes}, defined as those approaching
the center subspace of $\mA_\pm$ as $\lambda\to 0$,
are $v_f^\pm$, $v_r^\pm$ (resp.  $v_{f,j}^\pm$, $v_{r,i}^\pm$).
As {\it fast} growing and decaying modes, defined as those
approaching the stable and unstable subspace of $\mA_\pm$,
hence spectrally separated both from each other and from slow modes,
may always be chosen analytically in a neighborhood of $\lambda=0$
by the same argument used in Lemma \ref{Kat}, we obtain by our
asymptotic description of slow modes the following important extension.

\begin{lem}\label{enormal}
For Majda's model (more generally, the vectorial version including
reactive Navier--Stokes equations with artificial viscosity),
normal modes extend to $\Lambda\cup B(0,r)$
for $r>0$ sufficiently small,
with low-frequency (slow mode) asymptotics as described in Section \ref{cc}
(resp. \ref{syseval}).
\end{lem}

\begin{proof}
The explicit Taylor expansions of Section \ref{cc} (resp. \ref{syseval})
yield analytic extensions on $B(0,r)$ of slow modes $V_j^\pm$, spanning
invariant subspaces of $\mA_\pm$.
As fast modes always have such analytic extensions, we may combine them
to obtain analytic bases $V_j^\pm$ for invariant subspaces $S^\pm$, $U^\pm$
of $\mA_\pm$ extending those of Lemma \ref{Kat}.
These in turn determine $\mA_\pm$-invariant projections onto those
subspaces, which must therefore be the unique analytic extension of
the corresponding projections on $\Lambda \cap B(0,r)$, and thus
an analytic extension onto $\Lambda \cup B(0,r)$.  The result then
follows again by the result of Kato as in the proof of Lemma \ref{Kat}.
\end{proof}

\br\label{extensionrem}
\textup{
The extension of normal modes through the essential spectrum boundary
to a neighborhood of $\lambda=0$ is crucial for all that follows;
see \cite{GZ, ZH} for further discussion.
}
\er
\subsection{Construction of the Evans function}\label{evans}

\begin{defi}\label{d:evans}
\textup{
On any simply connected subset of the domain of consistent splitting,
let $V_{1}^+, \dots, V_{k}^+$ and
$V_{k+1}^-, \dots, V_{N}^-$
be analytic bases for $S_+$ and $U_-$, as described in Lemma \ref{Kat}.
Then, the {\it Evans function} for \eqref{gfirstorder}
associated with this choice of limiting bases is defined as
\ba \label{e:evans}
D(\lambda)&:=
\det
\Big(
 W_{1}^+, \dots,  W_k^+,
 W_{k+1}^-, \dots, W_{N}^-
\Big)_{|x=0, \lambda} \\
&=
\det
\Big(
P_+ V_{1}^+, \dots, P_+ V_k^+,
P_- V_{k+1}^-, \dots, P_- V_{N}^-
\Big)_{|x=0, \lambda},
\ea
where $P_\pm$ are the transformations described in Lemma \ref{conjugation}.
}
\end{defi}

\br\label{unique}
\textup{
Note that $D$ is independent of the choice of $P_\pm$; for,
by uniqueness of stable/unstable manifolds,
the exterior products (minors)
$P_+ V_{1}^+\wedge \dots\wedge P_+ V_k^+$ and
$P_- V_{k+1}^-\wedge \dots\wedge P_- V_{N}^-$
are uniquely determined by their behavior as $x\to + \infty$, $-\infty$,
respectively.
}
\er

\begin{prop}[MaZ3, Z]\label{evalue}
Both the Evans function and the stable/unstable
subspaces $\CalS^+$ and $\CalU^-$ are
analytic on the entire simply
connected subset of the domain of consistent splitting
on which they are defined.
Moreover, for $\lambda$ within this region,
equation \eqref{gfirstorder} admits
a nontrivial solution $W\in L^2(x)$ if and only if
$D(\lambda)=0$.
\end{prop}


\br\label{multiplicity}
\textup{
In the case that \eqref{gfirstorder} describes an eigenvalue
equation associated with an ordinary differential operator $L$,
$\lambda \in \CC^1$,
Proposition \ref{evalue} implies that eigenvalues of $L$ agree
in location with zeroes of $D$.
In \cite{GJ1, GJ2}, Gardner and Jones have shown that
they agree also in multiplicity; see also Lemma 6.1, \cite{ZH},
or Proposition 6.15 of \cite{MaZ3}.
}
\er

By Lemma \ref{enormal}, we have immediately that $D$ extends to
a neighborhood of the origin.

\begin{lem}\label{eD}
For Majda's model (more generally, the vectorial version including
reactive Navier--Stokes equations with artificial viscosity),
$D$ extends analytically to $\Lambda\cup B(0,r)$
for $r>0$ sufficiently small.
\end{lem}

\section{The resolvent kernel}\label{resolvent}
Next,
we estimate the resolvent kernel
$G_\lambda(x,y):=(L-\lambda)^{-1}\delta_y(x)$
associated with the linearized operator $L$ about the wave,
again following the abstract framework developed in \cite{ZH, MaZ3},
Rewriting \eqref{model} in vectorial form
\be\label{vnonlin}
U_t + F(U)_x+G(U)=BU_{xx},\\
\ee
\be\label{fluxes}
U:=(u, z), \qquad
F(U)=\bp f(u,z)-su\\ -s z \ep,
\qquad
G(U)= \bp \phi(u)qkz\\ -\phi(u)kz\ep,
\qquad
B=\bp b & 0 \\ 0 & d\ep,
\ee
in coordinates moving with a given traveling-wave
profile $\bar U(x)=(\bar u, \bar z)(x)$
(stationary, in the moving coordinate frame)
and linearizing about $\bar U$, we obtain the linearized equations
\be\label{vlin}
U_t=LU:= BU_{xx}- (AU)_x + CU,
\ee
where
$$
A:=dF(U)
=\bp \alpha & 0 \\
0 & -s
\ep,
\qquad
B=\bp b & 0 \\
0 & d
\ep,
\qquad
C:=dG(U)=\bp qk \phi'(\bar u)\bar z & qk\phi(\bar u) \\
-k \phi'(\bar u)\bar z & -k\phi(\bar u) \\
\ep.
$$

The results of \cite{ZH, MaZ3}
for general strictly parabolic systems of the form \eqref{vlin},
$U\in \RR^n$,
state that the resolvent kernel
is a meromorphic function on the domain of consistent splitting
defined in Section \ref{cc}, where it is determined by
\be\label{Glambda}
(L-\lambda) G_\lambda= \delta_y(x)
\ee
and the property of decay as $x$, $y\to \pm \infty$.
Moreover, they give an explicit description of $G_\lambda$
in terms of the normal modes of the eigenvalue equation
constructed in Section \ref{normal}, from which we may
extract sharp pointwise bounds.
We cite here the relevant theory, referring the reader to
\cite{ZH, MaZ3} for proof.

\subsection{Duality relation}
Consider solutions $U$ of the
eigenvalue equation $(L-\lambda)U=0$ and solutions $\tilde U$ of its adjoint
$(L^*-\bar \lambda)\tilde U$, where
\be\label{Ladj}
L^*\tilde U:= B^{\trans}\tilde U_{xx}+ A^\trans \tilde U_x + C^\trans \tilde U
\ee
denotes the $L^2$ adjoint of $L$
and $\bar \lambda$ the complex conjugate of $\lambda$.
Introducing the phase-variables $W:=(U,U')$ and
$\tilde W:=(\tilde U, \tilde U')$, write the eigenvalue
equation of $L$ and its adjoint as first-order ODE of form
\eqref{gfirstorder} in $W$ and $\tilde W$.
Then, we have the following key relation.

\begin{lem}[\cite{ZH, MaZ3}]\label{dual}
$W=(U,U')$ satisfies \eqref{gfirstorder} if and only if
\be\label{duality}
\tilde W^* \CalS W\equiv \text{\rm constant}
\ee
for all $\tilde W=(\tilde U,\tilde U')$ satisfying
the adjoint eigenvalue equation, and vice versa, where
\be\label{S}
\CalS:= \bp -A & B\\
-B & 0\ep,
\qquad
\CalS^{-1}=\
\bp 0 & -B^{-1} \\ B^{-1} & -B^{-1} A B^{-1} \ep.
\ee
\end{lem}

\begin{proof}
Property \eqref{duality} follows immediately from the relation
below,
 which we obtain from integration by parts:
$$
\tilde W^*\CalS W|^{x_2}_{x_1}= \langle(L-\lambda)^*\tilde U,
U\rangle_{L^2(x_1,x_2)} - \langle\tilde U,(L-\lambda)
U\rangle_{L^2(x_1,x_2)}=0,
$$
by the definition of the adjoint operator.
 (Indeed, the righthand
side may be viewed as defining the quadratic form $\CalS$.)
\end{proof}

\subsection{Domain of consistent splitting}\label{splitting}
Define
\be\label{Lambda}
\Lambda:=\cap \Lambda_j^\pm,\quad j=1,\dots,n,
\ee
where $\Lambda_j^\pm$ denote the open sets bounded on the left by
the algebraic curves $\lambda_j^\pm(\xi)$
determined by the eigenvalues of the symbols
$-\xi^2B_\pm - i\xi A_\pm + C_\pm$
of the limiting constant-coefficient operators
\be \label{limitL}
L_\pm U:=
B_\pm U'' -A_\pm U'+ C_\pm U
\ee
as $\xi$ is varied along the real axis.
The curves $\lambda_j^\pm(\cdot)$
comprise the essential spectrum of operators $L_\pm$.
For Majda's model,
the computations of Section \ref{cc} yield
\be \label{containment}
\Lambda\subset
\Omega_\eta:=
\{\lambda: \Re \lambda >
\max\{-\eta_1 |\Im \lambda|, \eta_2 |\Im \lambda|^2\} \},
\qquad
\eta_j>0.
\ee

\begin{lem}[\cite{MaZ3}]\label{frozen}
The set $\Lambda$ is equal to the
component containing real $+\infty$
of the domain of consistent splitting (defined in
Section \ref{cc}) for
the eigenvalue equation of $L$
written as a first-order ODE \eqref{gfirstorder}.
\end{lem}

\subsection{Basic solution formula}

Let
\be\label{phi+}
\Phi_j^+
:= P_+ V_j^+,\quad j=1\dots,n,
\ee
and
\be\label{phi-}
\Phi_{j}^- := P_- V_j^-, \quad j=n+1,\dots, 2n
\ee
denote the  locally analytic bases of the
stable manifold at $+\infty$
and the unstable manifold at $-\infty$ of solutions of
the eigenvalue equation \eqref{eval} written as a first-order system
\eqref{gfirstorder}
that were found in Section \ref{normal} (i.e., the normal modes), and set
\be\label{Phi+-}
\Phi^+:= (\Phi_1^+,\dots, \Phi_n^+),
\quad
\Phi^-:= (\Phi_{n+1}^-,\dots, \Phi_{2n}^-),
\ee
and
\be\label{Phi}
\Phi:= (\Phi^+,\Phi^-).
\ee
Define the solution operator from $y$ to $x$ of
\eqref{eval}, denoted by $\Fyx$, as
 \be\label{Fyx}
  \Fyx= \Phi(x,\lambda)\, \Phi^{-1}(y,\lambda)
 \ee
and the projections $\Pi^\pm_y$ on the stable manifolds at $\pm\infty$ as
  \be\label{Pi}
  \Pi^+_y=(\,\Phi^+(y,\lambda)\quad 0\,)\,\Phi^{-1}(y,\lambda)
   \quad\text{and}\quad
  \Pi^-_y=(\,0\quad \Phi^-(y,\lambda)\,)\,\Phi^{-1}(y,\lambda).
 \ee

Then, we have the following general result established in \cite{ZH, MaZ3}.

\begin{prop}[\cite{ZH, MaZ3}]\label{res1}
With respect to any $L^p$, $1\le p\le \infty$,
$\Lambda$ consists entirely of normal points of $L$, i.e.,
resolvent points, or isolated eigenvalues of constant multiplicity.
On this domain,
the resolvent kernel $G_\lambda$ is meromorphic,
with representation
  \be\label{basicflow}
  G_\lambda(x,y)=\begin{cases}
   (I_n,0)\Fyx\Pi^+_y\, \CalS^{-1}(y)
(I_n,  0)^{tr}
                                  \qquad &x>y,\\
   -(I_n,0)\Fyx\Pi^-_y\, \CalS^{-1}(y)
(I_n,  0)^{tr}
                                  \qquad &x<y,\\
                        \end{cases}
 \ee
$\CalS^{-1}$ as described in \eqref{S}.
Moreover, on any compact subset $K$ of $\rho(L)\cap \Lambda$
($\rho(L)$ denoting resolvent set), there hold the uniform decay estimates
\be\label{uKdecay}
|\partial_x^j\partial_y^k G_\lambda(x,y)|\le Ce^{-\eta |x-y|},
\ee
$0\le |j|+|k|\le 1$,
where $C>0$ and $\eta>0$ depend only on $K$, $L$.
\end{prop}

\br\label{augmented}
\textup{
Formula \eqref{basicflow} extends \cite{ZH, MaZ3}
to the full phase-variable representation
  \ba\label{augflow}
  &\bp
G_\lambda(x,y)&  (\partial/\partial y)G_\lambda(x,y)
\\
(\partial/\partial x)G_\lambda(x,y) &
(\partial/\partial x) (\partial/\partial y)G_\lambda(x,y)
\\
\ep
=
\begin{cases}
   \Fyx\Pi^+_y\, \CalS^{-1}(y)
                                  \qquad &x>y,\\
   -\Fyx\Pi^-_y\, \CalS^{-1}(y)
                                  \qquad &x<y.\\
                        \end{cases}
 \ea
}
\er

\subsection{Generalized spectral decomposition}
Formula \eqref{basicflow} suffices for the description
of intermediate- and high-frequency behavior.
For the treatment of the key low-frequency regime,
it is preferable to use a modified representation of the
resolvent kernel
consisting of a {\it scattering decomposition} in solutions
of the forward and adjoint eigenvalue equations.

{}From \eqref{duality}, it follows that if there are $n$ independent solutions
$\phi^+_1, \dots, \phi^+_{n}$ of $(L-\lambda I)U=0$ decaying at $+\infty$,
and $n$ independent solutions $\phi^-_{n+1}, \dots, \phi^-_{2n}$ of the same
equations decaying at $-\infty$, then there exist $n$ independent
solutions $\tilde\psi^+_{n+1}, \dots, \tilde\psi^+_{2n}$ of $(L^*-\lambda^* I){\tilde U}=0$
decaying at $+\infty$, and $n$ independent solutions
${\tilde\psi}^-_{1}, \dots, \tilde\psi^-_{n}$ decaying at $-\infty$.
Precisely, setting
\be\label{psi+}
\Psi_j^+
:= P_+ V_j^+,\quad j=n+1\dots,2n,
\ee
\be\label{psi-}
\Psi_{j}^- := P_- V_j^-, \quad j=1,\dots, n,
\ee
and
\be\label{Psi}
\Psi:= (\Psi^+,\Psi^-),
\ee
similarly as in \eqref{phi+}--\eqref{Phi},
where $\Psi_j^\pm=(\psi_j^\pm, (\psi_j^\pm)')$
are exponentially growing solutions obtained through
Lemma \ref{conjugation}, we may define dual exponentially decaying
and growing solutions $\tilde \psi_j^\pm$ and $\tilde \phi_j^\pm$ via
\be\label{dualbase}
\bp \tilde \Psi & \tilde \Phi \ep^* \CalS
\bp \Psi & \Phi \ep_\pm\equiv I.
\ee

\begin{cor}[\cite{ZKochel, MaZ3}]\label{kochelformulae}
On $\Lambda \cap \rho(L)$, there hold
\be\label{Mplusrep}
G_\lambda(x,y) =
\sum_{j,k} M^+_{jk}(\lambda ) \phi^+_j (x;\lambda) \Tpsi^-_k(y;\lambda)^*
\ee
for $y\le 0\le x$,
\be\label{dplusrep}
G_\lambda(x,y) =
\sum_{j,k} {d}^+_{jk}(\lambda ) \phi^-_j (x;\lambda)
\Tpsi^-_k(y;\lambda)^*
-\sum_k \psi^-_k (x;\lambda)
\Tpsi^-_k (y;\lambda)^*
\ee
for $y\le x\le 0$, and
\be\label{dminusrep}
G_\lambda(x,y) =
\sum_{j,k} {d}^-_{jk}(\lambda ) \phi^-_j (x;\lambda)
\Tpsi^-_k(y;\lambda)^* +
\sum_k \phi^-_k (x;\lambda)
\tilde \phi^-_k (y;\lambda)^*
\ee
for $x\le y\le 0$, with
\be\label{Mplus}
M^+=
(-I,0)
\bp
\Phi^+ & \Phi^-\\
\ep^{-1}
\Psi^-
\ee
and
\be\label{dplusminus}
d^\pm=
\big( 0 , I\big)
\bp
\Phi^+ & \Phi^-\\
\ep^{-1 }
\Psi^-.
\ee
Symmetric representations hold for $y\ge 0$.
\end{cor}

\begin{proof}
Rearrangement of \eqref{basicflow} using \eqref{dualbase};
see \cite{MaZ3}.
\end{proof}


\begin{rems}\label{r:cc}
\textup{
1.  This representation reflects
the classical duality principle (see, e.g. [ZH], Lemma 4.2)
that the transposition ${G}_{\lambda}^*(y,x)$ of the Green's function
$G_\lambda(x,y)$ associated with operator $(L-\lambda)$ should
be the Green's function for the adjoint operator $(L^*-\lambda^*)$.
}

\textup{
2. In the constant-coefficient case, with a choice of common bases
$\Psi^\pm=\Phi^\mp$ at $\pm \infty$,
\eqref{Mplusrep}--\eqref{dplusminus} reduce to the simple formula
  \be\label{cc2}
  G_\lambda(x,y)=\begin{cases}
   -\sum_{j=k+1}^{N} \phi_j^+(x;\lambda)\tilde \phi_j^{+*}(y;\lambda)
\qquad &x>y,\\
   \sum_{j=1}^k \phi_j^-(x;\lambda)\tilde \phi_j^{-*}(y;\lambda)
    \qquad &x<y,\\
                 \end{cases}
 \ee
where, generically,
$\phi_j^\pm$, $\tilde {\phi}_j^\pm$ may be taken as pure exponentials
\be\label{spectral}
\phi_j^\pm(x)\tilde \phi_j^{\pm*}(y)=e^{\mu_j^\pm(\lambda)(x-y)}
V_j^\pm(\lambda) \tilde{V}_j^{\pm*}(\lambda).
\ee
This reveals an analogy to the usual representation
obtained by Fourier transform solution.
We see that \eqref{dplusrep} (resp. \eqref{dminusrep}) in
the far-field limit $x\to +\infty$ (resp. $x\to -\infty$)
consists of the limiting constant-coefficient resolvent kernel plus
an exponentially decaying error term.
}
\end{rems}

\subsection{Resolvent kernel bounds}\label{resbounds}

The basic bound \eqref{uKdecay} is sufficient to treat intermediate
frequencies $r\le |\lambda|\le R$, $r$, $R>0$ (indeed, it is difficult
to say more in this regime).
As $|\lambda|\to \infty$, or $\lambda\to 0 \in \partial\rho(L)$,
however, the bounds are not uniform, and so separate
analyses are needed in these high- and low-frequency regimes.

In the high-frequency ($\sim$ short-time) regime,
we have the following classical
analytic semigroup-type bounds following from strict parabolicity alone.
These may be obtained from the basic formula
\eqref{basicflow} using the parabolic rescaling $x\to x|\lambda|^{1/2}$
and estimates (the Tracking Lemma of \cite{ZKochel, Z, MaZ3})
for slowly-varying-coefficient ODE;
see, e.g., \cite{Sa, ZH, ZKochel}.

\begin{prop} [High-frequency bounds \cite{ZH}]\label{HF}
For general strictly parabolic systems of form \eqref{vlin},
$R$, $C>0$ sufficiently large, and $\eta_1$, $\eta_2$, $\theta >0$
sufficiently small,
\ba\label{e:HF}
|G_{\lambda} (x,y)| &\le C |\lambda|^{-1/2}
 e^{-\theta |\lambda|^{\frac{1}{2}} |x-y|},
\qquad
|\partial_x G_{\lambda} (x,y)|,  \,
|\partial_y G_{\lambda} (x,y)| &\le C e^{-\theta
|\lambda|^{\frac{1}{2}}|x-y|} \\
\ea
for all $\lambda\in \Omega_\eta \setminus B(0,R)$, with
$\Omega_\eta$ as defined in \eqref{containment}.
\end{prop}

Thus, the only resolvent bounds
that depend on the details of the model are the
crucial low-frequency ($\sim$ large-time) bounds
carrying information relevant to large-time asymptotics.

For the Majda model, these are as follows.

\begin{prop}[Low-frequency bounds]\label{LF}
Let $\bar U=(\bar u, \bar z)$ be a traveling wave profile of
Majda's model, satisfying (\cD).
Then, for $r>0$ sufficiently small, the resolvent kernel $G_\lambda$ has
a meromorphic extension onto $B(0,r)\subset \CC$, which may be decomposed as
\be\label{ldecomp}
G_\lambda=E_\lambda + S_\lambda +R_\lambda,
\ee
where
\be \label{Elambdaminus}
E_\lambda(x,y)
:=
\begin{cases}
\lambda^{-1}
\bar U'(x) \pi_f^-(y)^\trans
e^{(\lambda/\alpha^{-} - \lambda^2 /{\alpha^{-}}^3 )y}
& \alpha^- > 0,\\
\lambda^{-1} \bar U'(x) \pi_r^-(y)^\trans
& \alpha^- < 0,\\
\end{cases}
\ee
for $y\le 0$,
with $\pi_j^-$ bounded solutions of the adjoint eigenvalue equation for
$\lambda=0$, $\pi_f^-$ convergent as $y\to - \infty$ to $cL_f^-$ and
$\pi_r^-$ exponentially decaying as $y\to -\infty$,
\be \label{Elambdaplus}
E_\lambda(x,y)
:=
\begin{cases}
\begin{aligned}
&\lambda^{-1}
\bar U'(x) \pi_f^+(y)^\trans
e^{(\lambda/\alpha^{+} - \lambda^2 /{\alpha^{+}}^3 )y}\\
&\quad  +\lambda^{-1}
\bar U'(x) \pi_r^+(y)^\trans
e^{(-\lambda/s + \lambda^2 d/s^3 )y}\\
\end{aligned}
& \alpha^+ < 0,\\
{}&{}\\
\lambda^{-1}
\bar U'(x) \pi_r^+(y)^\trans
e^{(-\lambda/s + \lambda^2 d/s^3 )y}
& \alpha^+ > 0,\\
\end{cases}
\ee
for $y\ge 0$,
with $\pi_j^+$ bounded solutions of the adjoint eigenvalue equation for
$\lambda=0$, convergent as $y\to + \infty$;
\be\label{tildeSlambda1}
S_\lambda(x,y):=
\begin{cases}
cR_f^{+}  {L_f^{-}}^t
e^{(-\lambda/\alpha^{+} + \lambda^2 /{\alpha^{+}}^3 )x
+(\lambda/\alpha^{-} - \lambda^2 /{\alpha^{-}}^3 )y}&
\alpha^-, \alpha^+>0\\
0 & \hbox{\rm otherwise}
\end{cases}
\ee
for $y\le 0\le x$, $R_f^-$ and $L_f^-$ constant vectors
as defined in \eqref{dualode},
\be\label{tildeSlambda2}
S_\lambda(x,y):=
\begin{cases}
R_f^{-}  {L_f^{-}}^\trans
e^{ (-\lambda/\alpha^{-} + \lambda^2 /{\alpha^{-}}^3 )(x-y)}
& \alpha^->0 \\
0 & \hbox{\rm otherwise}
\end{cases}
\ee
for $y\le x\le 0$,
and
\be \label{tildeSlambda3}
S_\lambda(x,y):=
\begin{cases}
R_f^{-}  {L_f^{-}}^\trans
e^{ (-\lambda/\alpha^{-} + \lambda^2 /{\alpha^{-}}^3 )(x-y)}
& \alpha^-<0 \\
0 & \hbox{\rm otherwise}
\end{cases}
\ee
for $ x\le y \le 0$, with similar relations for $y\ge 0$;
and $R_\lambda$ denotes a faster-decaying
residual term.\footnote{See, e.g., \cite{MaZ3, Z}
for bounds in the viscous Lax shock $\sim$ strong detonation case.}
\end{prop}

\begin{rems}\label{LFrems}
\textup{
1. Recall, \eqref{eq:v_r_taylor}--\eqref{fluidform}, \eqref{dualode}--\eqref{dfluidform},
that vectors $R_f^\pm=(*,0)^\trans$ appearing in scattering
terms $S_\lambda$ have vanishing $z$-component, and also
$L_F\perp(-q,1)$, a fact that will be
important in our later nonlinear stability analysis.
}

\textup{
2. The case $\alpha^-<0$ in \eqref{Elambdaminus} occurs only in the
extreme situation of a strong deflagration,
for which there are no incoming characteristics on the lefthand side
$y\le 0$.  This is essentially the only difference from the corresponding
proposition for viscous shock waves in the general Lax or
undercompressive case, and represents just an anomaly in bookkeeping.
}
\end{rems}

\begin{proof}
Similarly as in the viscous shock case treated in \cite{MaZ3}, this
follows from representations \eqref{Mplusrep}--\eqref{dplusminus}
of Corollary \ref{kochelformulae},
estimating modes $\phi_j^\pm$, $\psi_j^\pm$,
$\tilde \phi_j^\pm$, and $\tilde \psi_j^\pm$
using the asymptotic description given by the Conjugation Lemma
together with the constant-coefficient analysis of Section \ref{cc},
and estimating scattering coefficients $M_{jk}$, $d_{jk}^\pm$
crudely by Laurent series: e.g.,
$$
d_{jk}=d_{jk}^{-1}\lambda^{-1} +
d_{jk}=d_{jk}^{0}+ \dots,
$$
noting that pole terms of order $\lambda^{-k}$
correspond to zeroes of order $k$ of the Evans function, hence
(by (\cD)) are at most order $k=1$ and (without loss
of generality coordinatizing so that $\phi_1^+=\phi_{2n}^-=\bar U'(x)$)
involve only zero-eigenfunction $\bar U'(x)$ as $x$-dependent factor.
Specifically, $E_\lambda$ comprises exact pole terms, while $S_\lambda$
comprises order one terms involving products of slowly decaying
forward and dual modes (i.e., modes that are merely bounded for
$\lambda=0$), the latter estimated to exponentially decaying
error via the Conjugation Lemma, while $R_\lambda$ comprises remaining,
residual terms.

Vectors $R_j$, $L_j$ in the formulae for $S_\lambda$ arise through
the limiting, constant-coefficient analysis of Section \ref{cc}.
Finally, the information that $\pi_f^-$
at $\lambda=0$ converges as $y\to - \infty$ to $cL_f^-$ if $\alpha^->0$ and
to zero if $\alpha^-$ follows by the fact that in the first case
there exists but a single bounded, nondecaying solution of the
eigenvalue equation as $y\to -\infty$, with asymptotic direction $L_f^-$,
and in the second case there exists no bounded, nondecaying solution.
These facts, in turn, are readily verified by the constant-coefficient
analysis of Section \ref{cc}, combined with the Conjugation Lemma.
\end{proof}

\subsection{The system case}\label{sysres}
Proposition \ref{resbounds} admits a straightforward generalization
to the system case (i.e., the vectorial version of
\eqref{model} discussed in Remark \ref{modeling}.3),
substituting in place of $\alpha^\pm$, $L_f^\pm$, $R_f^\pm$
the eigenvalues $\alpha_j^\pm$ and eigenvectors $L_{f,j}^\pm$,
$R_{f,j}^\pm$ of $f_u(U_\pm)-sI$ and in place of $L_r^+$, $R_r^+$
the eigenvectors $L_{r,i}^+$, $R_{r,i}^+$ of the (now matrix-valued)
diffusion coefficient $d$, as described in Section \ref{syseval}
(recall, there are no slow reactive modes
on the minus infinity side $x\le 0$).
See Proposition 4.22, \cite{Z} for a corresponding description in the
viscous shock case.

\begin{prop}[Low-frequency bounds]\label{sysLF}
Let $\bar U=(\bar u, \bar z)$ be a traveling wave profile of
vectorial Majda's model under dissipativity hypothesis
\eqref{diss}, satisfying (\cD).\footnote{Recall, Remark \eqref{E-condition},
this implies in part that $\bar U$ is a unique, transversal connection.}
Then, for $r>0$ sufficiently small, the resolvent kernel $G_\lambda$ has
a meromorphic extension onto $B(0,r)\subset \CC$, which may be decomposed as
\be\label{sysldecomp}
G_\lambda=E_\lambda + S_\lambda +R_\lambda,
\ee
where
\be \label{sysElambdaminus}
E_\lambda(x,y)
:=
\begin{cases}
\lambda^{-1}
\sum_{\alpha_{j}^->0}\bar U'(x) \pi_{f,j}^-(y)^\trans
e^{(\lambda/\alpha_j^{-} - \lambda^2 b_j^- /{\alpha_j^{-}}^3 )y}
& \hbox{\rm some } \alpha_j^- > 0,\\
\lambda^{-1} \bar U'(x) \pi_r^-(y)^\trans
& \hbox{\rm all } \alpha_j^- < 0,\\
\end{cases}
\ee
for $y\le 0$,
with $\pi_{f,j}^-$ bounded solutions of the adjoint eigenvalue equation for
$\lambda=0$, exponentially convergent as $y\to - \infty$ to $c_j^-L_{f,j}^-$
and $\pi_{r,i}^-$ exponentially decaying as $y\to -\infty$,
\be \label{sysElambdaplus}
E_\lambda(x,y)
:=
\begin{cases}
\begin{aligned}
&\lambda^{-1}
\sum_{\alpha_{j}^+<0}
\bar U'(x) \pi_{f,j}^+(y)^\trans
e^{(\lambda/\alpha_j^{+} - \lambda^2 b_j^+ /{\alpha_j^{+}}^3 )y}\\
&\quad  +\lambda^{-1}
\sum_{i=1}^m
\bar U'(x) \pi_{r,i}^+(y)^\trans
e^{(-\lambda/s + \lambda^2 d_i^+/s^3 )y}\\
\end{aligned}
& \hbox{\rm some } \alpha_j^+ < 0,\\
{}&{}\\
\begin{aligned}
\lambda^{-1}
\sum_{i=1}^m
\bar U'(x) \pi_{r,i}^+(y)^\trans
e^{(-\lambda/s + \lambda^2 d_i^+/s^3 )y}
\end{aligned}
& \hbox{\rm all } \alpha_j^+ > 0,\\
\end{cases}
\ee
for $y\ge 0$,
with $\pi_{f,j}^+$, $\pi_r^+$ bounded solutions of the adjoint eigenvalue
equation for $\lambda=0$, exponentially convergent as $y\to + \infty$;
\be\label{systildeSlambda1}
S_\lambda(x,y):=
\sum_{\alpha_{k}^-, \alpha_{j}^+>0}
c^{j,+}_{k,-}
R_{f,j}^{+}  {L_{f,k}^{-}}^t
e^{(-\lambda/\alpha_j^{+} + \lambda^2 b_j^+ /{\alpha_j^{+}}^3 )x
+(\lambda/\alpha_j^{-} - \lambda^2 b_k^- /{\alpha_j^{-}}^3 )y}
\ee
for $y\le 0\le x$,
$R_{f,j}^-$ and $L_{f,k}^-$ constant vectors as defined in Section
\ref{syseval},
\ba\label{systildeSlambda2}
S_\lambda(x,y)&:=
\sum_{a_k^{-}>0}R_k^{-}  {L_k^{-}}^t
e^{ (-\lambda/\alpha^{-}_k + \lambda^2 b^{-}_k/{\alpha^{-}_k}^3 )(x-y)}\\
&\quad +
\sum_{\alpha_{f,j}^-<0, \alpha_k^->0}
c^{j,-}_{k,-}
R_{f,j}^{-}  {L_{f,k}^{-}}^\trans
e^{(-\lambda/\alpha^{-}_j + \lambda^2 b^{-}_j/{\alpha^{-}_j}^3 )x
+(\lambda/\alpha^{-}_k - \lambda^2 b^{-}_k/{\alpha^{-}_k}^3 )y}
\\
\ea
for $y\le x\le 0$,
and
\ba \label{systildeSlambda3}
S_\lambda(x,y)&:=
\sum_{\alpha_k^{-}<0}R_k^{-}  {L_k^{-}}^t
e^{ (-\lambda/\alpha^{-}_k + \lambda^2 b^{-}_k/{\alpha^{-}_k}^3 )(x-y)}
\\
&\quad +
\sum_{\alpha_{f,j}^-<0, \alpha_k^->0}
c^{j,-}_{k,-} R_{f,j}^{-}  {L_{f,k}^{-}}^\trans
e^{(-\lambda/\alpha^{-}_j + \lambda^2 b^{-}_j/{\alpha^{-}_j}^3 )x
+(\lambda/\alpha^{-}_k - \lambda^2 b^{-}_k/{\alpha^{-}_k}^3 )y}
\ea
for $ x\le y \le 0$,
where $c_j^\pm$, $c^{j,\pm}_{k,\pm}$ are scalar constants,
with similar relations for $y\ge 0$;
and $R_\lambda$ denotes a faster-decaying
residual term.
\end{prop}

\section{Green function bounds}\label{green}

We may now estimate the Green function
$G(x,t;y):=e^{Lt}\delta_y(x)$
associated with the linearized operator $L$ about the wave,
determined by
\be\label{G}
(\partial_t-L) G= 0,
\qquad
G(x,0;y)=\delta_y(x),
\ee
via the inverse Laplace-transform formula,
following the approach of \cite{ZH, MaZ3}.
We present our results using a bookkeeping scheme
similar to that of \cite{HZ} in the undercompressive
viscous shock case.

\subsection{Basic bounds}\label{bb}
Recall the standard notation
$$
\textrm{errfn} (z) := \frac{1}{2\pi} \int_{-\infty}^z e^{-\xi^2} d\xi.
$$

\begin{prop}\label{greenbounds}
Let $\bar U=(\bar u, \bar z)$ be a traveling wave profile of
Majda's model, satisfying (\cD).
Then, the Green function
$G(x,t;y)$ associated with the linearized equations \eqref{vlin}
may be decomposed as $G=E+\tilde G$, where
\begin{equation}\label{E}
E(x,t;y)= \bar U'(x) e(y,t),
\end{equation}
\begin{equation}\label{eminus}
  e(y,t)=
\begin{cases}
  \left(\textrm{errfn }\left(\frac{y+\alpha^{-}t}{\sqrt{4t}}\right)
  -\textrm{errfn }\left(\frac{y-\alpha^{-}t}{\sqrt{4t}}\right)\right)
  \pi_{f}^{-}(y)
& \alpha^- > 0,\\
 \pi_{r}^{-}(y)
& \alpha^- < 0,\\
\end{cases}
\end{equation}
for $y\le 0$ and
\be \label{eplus}
e(y,t)
:=
\begin{cases}
\begin{aligned}
&
  \left(\textrm{errfn }\left(\frac{y+\alpha^{+}t}{\sqrt{4t}}\right)
  -\textrm{errfn }\left(\frac{y+\alpha^{+}t}{\sqrt{4t}}\right)\right)
  \pi_{f}^{+}(y)
\\
&\quad  +
  \left(\textrm{errfn }\left(\frac{y+st}{\sqrt{4dt}}\right)
  -\textrm{errfn }\left(\frac{y-st}{\sqrt{4dt}}\right)\right)
  \pi_{r}^{+}(y)
\\
\end{aligned}
& \alpha_+ < 0,\\
{}&{}\\
  \left(\textrm{errfn }\left(\frac{y+st}{\sqrt{4dt}}\right)
  -\textrm{errfn }\left(\frac{y-st}{\sqrt{4dt}}\right)\right)
  \pi_{r}^{+}(y)
& \alpha^+ > 0,\\
\end{cases}
\ee
for $y\ge 0$, with $\pi_j^\pm$ as in Proposition \ref{LF}:
in particular,
\begin{equation}\label{pibounds}
|\pi_{j}^\pm|\le C, \qquad
|\partial_y \pi_{j}^\pm|\le Ce^{-\eta |y|},
\end{equation}
with $|\pi_f^-|\le Ce^{-\eta |y|}$ if $\alpha^-<0$ (strong deflagration case)
and $\pi_f^\pm\to L_f^\pm$ as $x\to -\infty$ otherwise,
and, denoting by $a_j^\pm$ the ``undamped'' characteristic speeds
$a_f^-=\alpha^-$, $a_f^+=\alpha^+$, and $a_r^+=-s$
(Note: $a_r^-$ does not appear),
\begin{equation}\label{Gbounds}
\begin{aligned}
|\partial_{x,y}^\alpha &\tilde G(x,t;y)|\le  Ce^{-\eta(|x-y|+t)}\\
& +\quad C(t^{-|\alpha|/2}+
|\alpha_x| e^{-\eta|x|}
+|\alpha_y| e^{-\eta|y|})\\
&\quad \times
\Big( \sum_{k}
t^{-1/2}e^{-(x-y-a_k^{-} t)^2/Mt} e^{-\eta x^+} \\
&+
\sum_{a_k^{-} > 0, \, a_j^{-} < 0}
\chi_{\{ |a_k^{-} t|\ge |y| \}}
t^{-1/2} e^{-(x-a_j^{-}(t-|y/a_k^{-}|))^2/Mt}
e^{-\eta x^+} \\
&+
\sum_{a_k^{-} > 0, \, a_j^{+}> 0}
\chi_{\{ |a_k^{-} t|\ge |y| \}}
t^{-1/2} e^{-(x-a_j^{+} (t-|y/a_k^{-}|))^2/Mt}
e^{-\eta x^-}
\\
&+
\sum_{a_k^{-} > 0}
\chi_{\{ |a_k^{-} t|\ge |y| \}}
t^{-1/2} e^{-(x+s(t-|y/a_k^{-}|))^2/Mt}
e^{-\eta |x|}\Big), \\
\end{aligned}
\end{equation}
$0\le |\alpha| \le 1$ for $y\le 0$ and symmetrically for $y\ge 0$,
for some $\eta$, $C$, $M>0$,
where $x^\pm$ denotes the positive/negative
part of $x$ and  indicator function $\chi_{\{ |a_k^{-}t|\ge |y| \}}$ is
$1$ for $|a_k^{-}t|\ge |y|$ and $0$ otherwise.
Moreover, for $x\le 0$, $|(0,1)\tilde G(x,t;y)|$
decays at the faster $x$-derivative rate $\alpha_x=1$,
as does $|(0,1)\tilde G(x,t;y) (1,0)^\trans|$ for any $x$,
and, for $y\le 0$,
$|\tilde G(x,t;y)(-q,1)^\trans |$
decays at the faster $y$-derivative rate $\alpha_y=1$.
\end{prop}

\begin{proof}
Reflecting the formal relation that
$G_\lambda$ is Laplace transform of $G$, we have the
Inverse Laplace transform formula
\be\label{ILT}
G(x,t;y)=
\frac{1}{2\pi i}
\oint_\Gamma e^{\lambda t} G_\lambda (x,y) \, d\lambda,
\ee
where $\Gamma= \partial \{ \lambda: \, Re \, \lambda > \theta_1- \theta_2
|\Im \, \lambda| \}$
is the boundary of  an appropriate sector containing the spectrum of $L$, $\theta_2>0$.
Following \cite{ZH, MaZ3},
we may thus convert the detailed resolvent kernel estimates of Proposition
\ref{resbounds} to estimates on the Green function via stationary phase,
or Riemann saddlepoint, estimates on \eqref{ILT}, exactly as was
done in the viscous shock case.

Specifically, using the property that $G_\lambda$ is meromorphic
on $\Omega_\eta \cup B(0,r)$ for $r$, $\eta>0$ sufficiently
small (see Propositions \ref{res1}, \ref{HF}, and \ref{LF}) and
analytic on the the resolvent set $\rho(L)$,
we may estimate the contribution of each of the various meromorphic
components $C_\lambda$ of $G_\lambda$ by a combination of direct evaluation
using Calculus of residues and strategic deformation of the contour $\Gamma$
so as to minimize
$$
\oint_\Gamma |e^{\lambda t} C_\lambda (x,y)| \, |d\lambda|=
\oint_\Gamma e^{\Re \lambda t}| C_\lambda (x,y)| \, |d\lambda|
$$
for each fixed $x$, $y$, $t$, with the main contribution to $E$ coming
from explicit evaluation of the corresponding low-frequency term $E_\lambda$
in Propositon \ref{LF} and the main contribution to $\tilde G$
coming from explicit evalution of $S_\lambda$.
See \cite{ZH, MaZ3, Z} for details.

It remains only to verify the key properties of faster decay of
$|(0,1)\tilde G(x,t;y)|$ for $x\le 0$,
$|(0,1)\tilde G(x,t;y) (1,0)^\trans|$ for general $x$,
and $|\tilde G(x,t;y)(-q,1)^\trans |$ for $y\le 0$.
To see the first property, we have only to observe,
in the bounds of Proposition \ref{resbounds},
that, for $x\le 0$, only fluid modes appear
in the rate-determining term $S_\lambda$, and
these lie in direction $R_f^-=(*,0)^\trans$ having
vanishing $z$-component.\footnote{The restriction
$x\le 0$ is necessary because of incoming (i.e.,
leftmoving) undamped
reaction waves for the case $y\ge 0$ not listed
in Proposition \ref{greenbounds}.}
For $x\ge 0$, fluid modes again lie in direction
$R_f^+=(*,0)^\trans$ having vanishing $z$-component,
while reactive terms appear as scalar multiples of
projector $R_r^+(L_r^+)^\trans=(0,*)$ orthogonal to
$(1,0)^\trans$.  Thus,
$|(0,1)\tilde G(x,t;y) (1,0)^\trans|$ to lowest order
involves only $z$-components of fluid terms, hence again
is faster decaying;
this yields the second property.
Likewise, for $y\le 0$, $S_\lambda (-q,1)^\trans=0$,
since $L_f^-\perp (-q,1)^\trans$, and this yields
the third property, completing the proof.
\end{proof}

\begin{rems}\label{outgoing}
\textup{
1. Similarly as in the viscous or relaxation shock case, the bounds
of Proposition \ref{greenbounds} may be interpreted as describing
the evlution of an initial delta-function perturbation at $y$
as the superposition of signals convecting along hyperbolic characteristics
and diffusing as approximate Gaussians until they strike the shock layer,
whereupon they scatter as reflected and transmitted waves along outgoing
characteristics, at the same time exciting the stationary mode $\bar U'$.
The main new feature in the reacting as compared to the nonreacting
case is the exiting signal along the reaction characteristic on the
lefthand ($x\le 0$) side,
for which a constant-coefficient analysis indicates that the Gaussian
signal is now exponentially decaying in time, due to burning of the
reactant.
This is reflected in the final term of \eqref{Gbounds}, consisting
of an ordinary Gaussian reflected left into an exponentially
penalized field $e^{-\theta|x|}$, a term indistinguishable in modulus
bound from a Gaussian multiplied by a factor decaying
exponentially in the travel time after reflection.
We call this leftgoing reactive characteristic speed ``damped''
and all others ``undamped'' to distinguish this behavior.
}

\textup{
2. A second difference between the reacting and nonreacting
case, this time confined to (scalar) Majda's model,
is the different structure of the excited term $E$ in the case
$\alpha^-<0$, a book-keeping anomaly arising because of the
absence of incoming waves on the lefthand side in this
(strong deflagration) case.
This different structure has essentially no effect on the analysis;
see Remark \ref{sdeflcase}.
}

\textup{
3. The improved bounds for $|(0,1)\tilde G(x,t;y)|$
and $|\tilde G(x,t;y)(-q,1)^\trans|$ for $x\le 0$
are similar to those of the relaxation case \cite{MaZ1},
to which the linearized equations are analogous on the
side $x\le 0$ on which $\phi>0$, and play a similarly important role
in the later nonlinear stability analysis.
In conservative coordinates $w:=u+qz$, $z$
of \eqref{eq:cmm1}--\eqref{eq:cmm2},
these bounds have the simpler statement that the Green function decays
more rapidly in its $z$-components, both output and input.
}
\end{rems}

\subsection{The system case}\label{sysbb}

The somewhat cumbersome summation notation of Proposition \ref{greenbounds}
is designed for easy generalization to the system case.  Indeed,
starting from Remark \ref{LFrems}.3, it is straightforward
to verify the analogous theorem for the full reactive
Navier--Stokes equations with artificial viscosity-- more generally,
the abstract vectorial model $u\in \RR^n$, $z\in \RR^m$
described in Remark \ref{modeling}.3-- with
undamped characteristic modes now
$a_j^\pm= \alpha_{1}^\pm,\dots, \alpha_{n}^\pm$, $-s$
($-s$ multiplicity $m=\dim z$), $\alpha_j^\pm$ and
$L_f^\pm$ denoting the eigenvalues and left eigenvectors of
$(\partial f/\partial u)(U_\pm)-sI$.
See \cite{HZ} for the analogous description in the viscous shock case.

\begin{prop}\label{sysgreenbounds}
Let $\bar U=(\bar u, \bar z)$ be a traveling wave profile of
vectorial Majda's model under dissipativity hypothesis
\eqref{diss}, satisfying (\cD).\footnote{Recall, \eqref{E-condition},
this implies in part that $\bar U$ is a unique, transversal connection.}
Then, the Green function
$G(x,t;y)$ associated with the linearized equations \eqref{vlin}
may be decomposed as $G=E+\tilde G$, where
\begin{equation}\label{sysE}
E(x,t;y)= \bar U'(x) e(y,t),
\end{equation}
\begin{equation}\label{syseminus}
  e(y,t)=
\begin{cases}
\sum_{\alpha_{j}^->0}
  \left(\textrm{errfn }\left(\frac{y+\alpha_j^{-}t}{\sqrt{4b_j^-t}}\right)
  -\textrm{errfn }\left(\frac{y-\alpha_j^{-}t}{\sqrt{4b_j^-t}}\right)\right)
 \pi_{f,j}^-(y)^\trans
& \hbox{\rm some } \alpha_j^- > 0,\\
 \pi_{r}^{-}(y)
& \hbox{\rm all } \alpha_j^- < 0,\\
\end{cases}
\end{equation}
for $y\le 0$ and
\be \label{syseplus}
e(y,t)
:=
\begin{cases}
\begin{aligned}
&
\sum_{\alpha_{j}^+<0}
  \left(\textrm{errfn }\left(\frac{y+\alpha_j^{+}t}{\sqrt{4b_j^+t}}\right)
  -\textrm{errfn }\left(\frac{y+\alpha_j^{+}t}{\sqrt{4b_j^+t}}\right)\right)
\pi_{f,j}^+(y)^\trans
\\
&\quad  +
\sum_{i=1}^m
  \left(\textrm{errfn }\left(\frac{y+st}{\sqrt{4d_i^+t}}\right)
  -\textrm{errfn }\left(\frac{y-st}{\sqrt{4d_i^+t}}\right)\right)
\pi_{r,i}^+(y)^\trans
\\
\end{aligned}
& \hbox{\rm some } \alpha_j^+ < 0,\\
{}&{}\\
\sum_{i=1}^m
  \left(\textrm{errfn }\left(\frac{y+st}{\sqrt{4d_i^+t}}\right)
  -\textrm{errfn }\left(\frac{y-st}{\sqrt{4d_i^+t}}\right)\right)
\pi_{r,i}^+(y)^\trans
& \hbox{\rm all } \alpha_j^+ > 0,\\
\end{cases}
\ee
for $y\ge 0$, with $\pi_j^\pm$ as in Proposition \ref{LF}:
in particular,
\begin{equation}\label{syspibounds}
|\pi_{j}^\pm|\le C, \qquad
|\partial_y \pi_{j}^\pm|\le Ce^{-\eta |y|},
\end{equation}
with $|\pi_f^-|\le Ce^{-\eta |y|}$ if $\alpha^-<0$ (strong deflagration case)
and $\pi_f^\pm\to L_f^\pm$ as $x\to -\infty$ otherwise,
and $\tilde G$ satisfies \eqref{Gbounds}.
Moreover, for $x\le 0$, $|(0,I_r)\tilde G(x,t;y)|$
decays at the faster $x$-derivative rate $\alpha_x=1$,
as does $|(0,I_r)\tilde G(x,t;y) (I_n,0)^\trans|$ for any $x$,
and, for $y\le 0$,
$|\tilde G(x,t;y)(-q^\trans,I_r)^\trans |$
decays at the faster $y$-derivative rate $\alpha_y=1$.
\end{prop}

\subsection{Linearized stability criterion}\label{linstab}


\begin{proof}[Proof of Theorem \ref{D}]
Sufficiency of (\cD) for linearized orbital stability follows
immediately by the bounds of Theorem \ref{greenbounds} (resp.
Remark \ref{outgoing}.4)
and standard $L^q\to L^p$ convolution bounds,
exactly as in the viscous shock case, setting
$$
\delta(t):= \int_{-\infty}^{+\infty}E(x,t;y)u_0(y)dy
$$
so that
$$
U-\delta(t)\bar U'= \int_{-\infty}^{+\infty}\tilde G(x,t;y)u_0(y)dy;
$$
see \cite{ZH, MaZ3, Z} for further details.
Necessity follows from more general spectral considerations not requiring
the detailed bounds of Theorem \ref{greenbounds};
see the discussion of effective spectrum in \cite{ZH, MaZ3, Z}.
The argument goes again exactly as in the viscous shock case.
\end{proof}

\section{Nonlinear stability}\label{nstab}
We can now readily establish nonlinear stability
by a combination of the methods used in \cite{HZ} to treat
general undercompressive viscous shock waves ($\sim$ $x\ge 0$ behavior),
and the methods used in \cite{MaZ1} to treat relaxation shocks
($\sim$ $x\le 0$ behavior).
As it costs no additional effort in bookkeeping,
we carry out this part of the argument in the full generality
of the system case.
Recall (Remark \ref{modeling}.3) that this includes
the artificial viscosity version of
the full reactive Navier--Stokes equations with multi-species
reaction and reaction-dependent equation of state.

Denoting by $a_j^\pm$ the ``undamped'' characteristic speeds:
$a_f^-=\alpha^-$, $a_f^+=\alpha^+$, and $a_r^+=-s$ for
Majda's model;
$a_{f,i}^-=\alpha_i^-$, $i=1,\dots, n$,
$a_{f,i}^+=\alpha_i^+$, $i=1,\dots, n$,
$a_{r,i}^+=-s$, $i=1,\dots, m$ in the system case $u\in \RR^n$, $z\in \RR^m$;
define
\begin{equation}\label{theta}
\theta(x,t):=
\sum_{a_j^-<0}(1+t)^{-1/2}e^{-|x-a_j^-t|^2/Lt}
+ \sum_{a_j^+>0}(1+t)^{-1/2}e^{-|x-a_j^+t|^2/Lt},
\end{equation}
\begin{equation}\label{psi1}
\begin{aligned}
\psi_1(x,t)&:=
\chi(x,t)\sum_{a_j^-<0}
(1+|x|+t)^{-1/2} (1+|x-a_j^-t|)^{-1/2}\\
&\quad+
\chi(x,t)\sum_{a_j^+>0}
(1+|x|+t)^{-1/2} (1+|x-a_j^+t|)^{-1/2},\\
\end{aligned}
\end{equation}
and
\begin{equation}\label{psi2}
\begin{aligned}
\psi_2(x,t)&:=
(1-\chi(x,t)) (1+|x-a_1^-t|+t^{1/2})^{-3/2}\\
&\quad +(1-\chi(x,t)) (1+|x-a_n^+t|+t^{1/2})^{-3/2},
\end{aligned}
\end{equation}
where $L>0$ is a sufficiently large constant and
$\chi(x,t)=1$ for
\be\label{chi}
x\in [\min_j \{a_j^- t, 0 \}, \max_j \{a_j^+ t, 0 \} ],
\ee
that is, for $x$ between the extremal outgoing undamped characteristics,
and zero otherwise.\footnote{This repairs a minor omission in \cite{HZ},
where \eqref{chi} was stated incorrectly as
$x\in [\min_j \{a_j^- t \}, \max_j \{a_j^+ t \} ]$.
The formulae differ in the case that there are no outgoing
characteristics on one side: extreme Lax shock or strong detonation.}

Then, we have the following pointwise version of Theorem \ref{nonlin}.

\begin{prop}\label{pwnonlin}
Let $\bar U(x-st)$ be a traveling combustion wave
of Majda's model (more generally, the vectorial version including
reactive Navier--Stokes equations with artificial viscosity)
and $|U_0(x)|\le E_0 (1+|x|)^{-3/2}$, $E_0$ sufficiently small.
Then, there exist $\delta(\cdot)$ and $\delta(+\infty)$ such that
\begin{equation}
\label{pointwise}
\begin{aligned}
|\tilde U(x,t)-\bar U^{\delta(t)}(x)|&\le C E_0
(\theta+\psi_1+\psi_2)(x,t),\\
  |\dot \delta (t)|&\le C E_0 (1+t)^{-1},\\
  |\delta(t)-\delta(+\infty)|&\le C E_0(1+t)^{-1/2},
\end{aligned}
\end{equation}
where $\tilde U$ denotes the solution of the same equations
with perturbed initial data $\tilde U_0=\bar U+U_0$.
\end{prop}

As discussed in the introduction,
we establish Proposition \ref{pwnonlin} by a combination
of the analysis of undercompressive viscous shock waves in \cite{HZ} and
of relaxation shocks in \cite{MaZ1}.

Following \cite{HZ}, set
\begin{equation}\label{pert}
U(x,t):=\tilde U(x+\delta(t), t)-\bar U(x),
\end{equation}
so that \eqref{vnonlin} becomes by Taylor expansion of
$F$, $G$:
\begin{equation}
\label{perteq}
U_t-LU=Q(U)_x+ R(U)+ \dot \delta (t)(\bar U_x + U_x),
\end{equation}
$L$ as in \eqref{vlin}, where
$R(U)=(-q,1)^\trans r(U)$, $r(U)$ scalar, with
\ba\label{Q}
Q(U)&=\mathcal{O}(|U|^2), \\
r(U)&=\mathcal{O}(|U|^2 e^{-\eta x^+}),\\
\ea
so long as $|U|$ remains bounded, where $x^+$ denotes the
positive part of $x$ and $\eta>0$.

\br\label{emphrem}
\textup{
Here, in the description of $R$, $r$ we have used the specific form
\be\label{rkphi}
G(U)=-\phi(u)\bp -qk \\ k \ep
\ee
of the reactive source in \eqref{vnonlin}, together with Taylor
expansion
\ba\label{sourcebound}
(\phi(\bar u+u)(\bar z + z)-
(\phi(\bar u)(\bar z)-
(\phi'(\bar u)u\bar z + \phi(\bar u)z)&=
\phi'(\bar u) uz + \phi''(\bar u+ \theta u)u^2\bar z,\\
\ea
$0< \beta <1$, and the fact that $\phi'(\bar u + v)\le Ce^{-\eta x^+}$
for $|v|$ sufficiently small, by assumption $u_+\not \in [u_i,u^i]$,
the property that $\phi'(u)\equiv 0$ for $u\not \in [u_i,u^i]$,
and exponential convergence of $\bar U(x)$ to $U_+$ as $x\to +\infty$.
This computation, and its exploitation in the later argument
(see especially the auxiliary bounds of Lemma \ref{auxconvolutions}),
are the main new features in the combustion context as compared
to the undercompressive viscous shock wave case.
}
\er

Recalling the standard fact that $\bar U'$ is a stationary
solution of the linearized equations \eqref{vlin},
$L\bar U'=0$, or
$$
\int^\infty_{-\infty}G(x,t;y)\bar U_x(y)dy=e^{Lt}\bar U_x(x)
=\bar U'(x),
$$
we have by Duhamel's principle:
$$
\begin{array}{l}
  \displaystyle{
  U(x,t)=\int^\infty_{-\infty}G(x,t;y)U_0(y)\,dy } \\
  \displaystyle{\qquad
  +\int^t_0 \int^\infty_{-\infty} G(x,t-s;y)(-q,1)^\trans
  r(U) (y,s)\,dy\,ds}\\
  \displaystyle{\qquad
  -\int^t_0 \int^\infty_{-\infty} G_y(x,t-s;y)
  (Q(U)+\dot \delta U ) (y,s)\,dy\,ds + \delta (t)\bar U'(x).}
\end{array}
$$
Defining
\begin{equation}
 \begin{array}{l}
  \displaystyle{
  \delta (t)=-\int^\infty_{-\infty}e(y,t) U_0(y)\,dy }\\
  \displaystyle{\qquad
  -\int^t_0 \int^\infty_{-\infty} e(y,t-s)(-q,1)^\trans
  r(U) (y,s)\,dy\,ds}\\
  \displaystyle{\qquad
  +\int^t_0\int^{+\infty}_{-\infty} e_{y}(y,t-s)(Q(U)+
  \dot \delta\, U)(y,s) dy ds, }
  \end{array}
 \label{delta}
\end{equation}
following \cite{Z3, MaZ1, MaZ2, MaZ4},
where $e$ is defined as in \eqref{eminus}--\eqref{eplus}
(that is, $e=\sum_j e_j$),
and recalling the decomposition $G=E+\tilde G$,
we obtain finally the {\it reduced equations}:
\begin{equation}
\begin{array}{l}
 \displaystyle{
  U(x,t)=\int^\infty_{-\infty} \tilde G(x,t;y)U_0(y)\,dy }\\
  \displaystyle{\qquad
  +\int^t_0 \int^\infty_{-\infty} \tilde G(x,t-s;y)(-q,1)^\trans
  r(U) (y,s)\,dy\,ds}\\
 \displaystyle{\qquad
  -\int^t_0\int^\infty_{-\infty}\tilde G_y(x,t-s;y)(Q(U)+
  \dot \delta U)(y,s) dy \, ds, }
\end{array}
\label{u}
\end{equation}
and, differentiating (\ref{delta}) with respect to $t$,
and observing that
$e_y (y,s)\rightharpoondown 0$ as $s \to 0$, as the difference of
approaching heat kernels:
\begin{equation}
 \begin{array}{l}
 \displaystyle{
  \dot \delta (t)=-\int^\infty_{-\infty}e_t(y,t) U_0(y)\,dy }\\
  \displaystyle{\qquad
  +\int^t_0 \int^\infty_{-\infty} e_t(y,t-s)(-q,1)^\trans
  r(U) (y,s)\,dy\,ds}\\
 \displaystyle{\qquad
  +\int^t_0\int^{+\infty}_{-\infty} e_{yt}(y,t-s)(Q(U)+
  \dot \delta U)(y,s)\,dy\,ds. }
 \end{array}
\label{deltadot}
\end{equation}
\medskip

The following integral estimates are established in \cite{HZ}.

\begin{lem}[Linear
estimates \cite{HZ}]\label{iniconvolutions}
Under the assumptions of Theorem \ref{nonlin},
\begin{equation}\label{iniconeq}
\begin{aligned}
\int_{-\infty}^{+\infty}|\tilde G(x,t;y)|(1+|y|)^{-3/2}\, dy
&\le C(\theta+\psi_1+\psi_2)(x,t),\\
\int_{-\infty}^{+\infty}|e_t(y,t)|(1+|y|)^{-3/2}\, dy
&\le C(1+t)^{-3/2},\\
\int_{-\infty}^{+\infty}|e(y,t)|(1+|y|)^{-3/2}\, dy
&\le C,\\
 \int^{+\infty}_{-\infty} |e(y,t)-e(y,+\infty)| (1+|y|)^{-3/2}\, dy
&\le C(1+t)^{-1/2},\\
\end{aligned}
\end{equation}
for $0\le t\le +\infty$,
some $C>0$, where $\tilde G$ and
$e$ are defined as in Proposition \ref{greenbounds}.
\end{lem}

\begin{lem}[Nonlinear
estimates \cite{HZ}]\label{convolutions}
Under the assumptions of Theorem \ref{nonlin},
\begin{equation}\label{coneq}
\begin{aligned}
\int_0^t\int_{-\infty}^{+\infty}|\tilde G_y(x,t-s;y)|\Psi(y,s)\, dy ds
&\le C(\theta+\psi_1+\psi_2)(x,t),\\
\int_0^t\int_{-\infty}^{+\infty}|e_{yt}(y,t-s)|\Psi(y,s)\, dy ds
&\le C(1+t)^{-1},\\
\int_t^{+\infty} \int_{-\infty}^{+\infty}|e_y(y,+\infty)|\Psi(y,s)\, dy
&\le C\gamma(1+t)^{-1/2},\\
\int_0^t\int_{-\infty}^{+\infty}
|e_y(y,t-s)- e_y(y,+\infty)| \Psi(y,s)\, dyds
&\le C(1+t)^{-1/2},\\
\end{aligned}
\end{equation}
for $0\le t\le +\infty$,
some $C>0$, where $\tilde G$ and
$e$ are defined as in Proposition \ref{greenbounds} and
\begin{equation}\label{source}
\begin{aligned}
\Psi(y,s)&:=
(1+s)^{1/2}s^{-1/2}(\theta + \psi_1+\psi_2)^2(y,s)\\
&\qquad +
(1+s)^{-1} (\theta+\psi_1+\psi_2)(y,s).
\end{aligned}
\end{equation}
\end{lem}

\br\label{sdeflcase}
\textup{
The case $\sigma(\alpha^-)<0$ in \eqref{eminus}, \eqref{syseminus},
occurring for
strong deflagrations as described in Remark \ref{LFrems}.2,
is the only one requiring discussion, since in all other cases
the bounds are identical to those of the shock case.
We have only to note that $|e|$, $|e_y|$, $|e_{yt}|=|e_t|\equiv 0$
in this case also satisfies the same bounds (or better) that are actually used
in the proofs of Lemmas \ref{iniconvolutions} and \ref{convolutions}.
}
\er

To these, we add the following auxiliary estimates special to
the combustion case.

\begin{lem}[Auxiliary estimates] \label{auxconvolutions}
Under the assumptions of Theorem \ref{nonlin},
\begin{equation}\label{auxconeq}
\begin{aligned}
\int_0^t\int_{-\infty}^{+\infty}|\tilde G(x,t-s;y)(-q,1)^\trans
e^{-\eta y^+}|\Psi(y,s)\, dy ds
&\le C(\theta+\psi_1+\psi_2)(x,t),\\
\int_0^t\int_{-\infty}^{+\infty}|e_{t}(y,t-s)(-q,1)^\trans
e^{-\eta y^+}|\Psi(y,s)\, dy ds
&\le C(1+t)^{-1},\\
\int_t^{+\infty} \int_{-\infty}^{+\infty}|e(y,+\infty)(-q,1)^\trans
e^{-\eta y^+}|\Psi(y,s)\, dy
&\le C\gamma(1+t)^{-1/2},\\
\int_0^t\int_{-\infty}^{+\infty}
|\big(e(y,t-s)- e(y,+\infty)\big)(-q,1)^\trans
e^{-\eta y^+}|\Psi(y,s)\, dy ds
&\le C(1+t)^{-1/2},\\
\end{aligned}
\end{equation}
for $0\le t\le +\infty$,
some $C>0$, $\tilde G$ and
$e$ as in Proposition \ref{greenbounds} and
$\Psi$ as in \eqref{source}.
\end{lem}

\begin{proof}
For the $\tilde G$-estimate, we have only to recall that,
by the bounds of Proposition \ref{greenbounds},
$\tilde G(-q, 1)^\trans$ obeys the bounds of $\tilde G_y$
for $y\le 0$, while $\tilde G e^{-\eta y^+}$ obeys the
bounds of $\tilde G_y$ for $y\ge 0$.
The $e$-estimates follow similarly, by the observation that,
for $y\le 0$, $e (-q,1)^\trans$
and $e_t (-q,1)^\trans$ decay like $(e^{-\eta|y|}+(1+t)^{-1/2})$
times the bounds for $e$ and $e_t$, hence obey the bounds for $|e_y|$
and $e_{yt}$, since $\pi_f^-$ is asymptotically
parallel to $L_f^-\perp (-q,1)^\trans$, with convergence at exponential
rate.
Likewise, for $y\ge 0$,
$e e^{-\eta y^+}$ and $e_t e^{-\eta y^+}$
decay like $(e^{-\eta|y|}+(1+t)^{-1/2})$ times the bounds
for $e$ and $e_t$, hence obey the bounds for $|e_y|$
anb $|e_{yt}$.
Thus, all bounds follow by the same arguments as in the proof of
Lemma \ref{convolutions}.
\end{proof}

\medskip

\begin{proof}[Proof of Proposition \ref{pwnonlin}]
With these observations, the proof of nonlinear stability
goes essentially as in \cite{HZ}.
Define
\begin{equation}
\label{zeta2}
 \zeta(t):= \sup_{y, 0\le s \le t}
 \Big( |U|(\theta+\psi_1+\psi_2)^{-1}(y,t)
 + |\dot \delta (s)|(1+s) \Big).
\end{equation}
We shall establish:

{\it Claim.} For all $t\ge 0$ for which a solution exists with
$\zeta$ uniformly bounded by some fixed, sufficiently small constant,
there holds
\begin{equation}
\label{claim}
\zeta(t) \leq C_2(E_0 + \zeta(t)^2).
\end{equation}
\medskip

{}From this result, provided $E_0 < 1/4C_2^2$,
we have that $\zeta(t)\le 2C_2E_0$ implies
$\zeta(t)< 2C_2E_0$, and so we may conclude
by continuous induction that
 \begin{equation}
 \label{bd}
  \zeta(t) < 2C_2E_0
 \end{equation}
for all $t\geq 0$.
(By standard short-time existence for artificial viscosity systems
(see, e.g., \cite{HoS, ZH}),
 $U\in C^1$ exists and $\zeta$ remains
continuous so long as $\zeta$ remains bounded by some uniform constant,
hence \eqref{bd} is an open condition.)
Thus, it remains only to establish the claim above.
\medskip

{\it Proof of Claim.}
We must show that $U(\theta+\psi_1+\psi_2)^{-1}$ and
$|\dot \delta(s)|(1+s)$
are each bounded by
$C(E_0 + \zeta(t)^2)$,
for some $C>0$, all $0\le s\le t$, so long as $\zeta$ remains
sufficiently small.

Recalling definition \eqref{zeta2},
we obtain for all $t\ge 0$ and some $C>0$ that
\begin{equation}\label{ubounds}
\begin{aligned}
|\dot \delta(t)|&\le \zeta(t)(1+t)^{-1},\\
|U(x,t)| &\le \zeta (t)(\theta +\psi_1+\psi_2)(x,t),\\
\end{aligned}
\end{equation}
and therefore
\begin{equation}\label{Nbounds}
\begin{aligned}
|(Q(U)+ \dot \delta U)(y,s)|&\le
C\zeta(t)^2 \Psi(y,s)
\end{aligned}
\end{equation}
with $\Psi$ as defined in \eqref{source}, for $0\le s\le t$.

Combining \eqref{Nbounds} with representations
(\ref{u})--(\ref{deltadot}) and
applying Lemmas \ref{iniconvolutions} and \ref{convolutions}, we obtain
$$
 \begin{aligned}
  |U(x,t)| &\le
  \int^\infty_{-\infty} |\tilde G(x,t;y)| |U_0(y)|\,dy
   \\
  &\qquad
  +\int^t_0 \int^\infty_{-\infty} |e(y,t-s)(-q,1)^\trans|
  |r(U) (y,s)|\,dy\,ds\\
 &\qquad +\int^t_0
  \int^\infty_{-\infty}|\tilde G_y(x,t-s;y)||(Q(U)+
  \dot \delta U)(y,s)| dy \, ds \\
  & \le
  E_0 \int^\infty_{-\infty} |\tilde G(x,t;y)|(1+|y|)^{-3/2}\,dy
   \\
  &\quad
C\zeta(t)^2 \int^t_0
   \int^\infty_{-\infty} |e(y,t-s)(-q,1)^\trans e^{-\eta y^+}|
  \Psi(y,s)|\,dy\,ds\\
 &\quad +
C\zeta(t)^2 \int^t_0
  \int^\infty_{-\infty}|\tilde G_y(x,t-s;y)|
\Psi(y,s) dy \, ds \\
&\le
C(E_0+\zeta(t)^2)(\theta + \psi_1+\psi_2)(x,t)
\end{aligned}
$$
and, similarly,
$$
\begin{aligned}
 |\dot \delta(t)| &\le \int^\infty_{-\infty}|e_t(y,t)|
  |U_0(y)|\,dy \\
&\qquad
  +\int^t_0 \int^\infty_{-\infty} |e_t(y,t-s)(-q,1)^\trans|
  |r(U) (y,s)|\,dy\,ds\\
  &\qquad +\int^t_0\int^{+\infty}_{-\infty} |e_{yt}(y,t-s)||(Q(U)+
 \dot \delta U)(y,s)|\,dy\,ds\\
&\le \int^\infty_{-\infty}E_0 |e_t(y,t)|(1+|y|)^{-3/2}\, dy\\
&\quad +
\int^t_0\int^{+\infty}_{-\infty}C\zeta(t)^2 |e_{t}(y,t-s)(-q,1)^\trans
e^{-\eta y^+}|\Psi(y,s) \,dy\,ds\\
\\
&\quad +
\int^t_0\int^{+\infty}_{-\infty}C\zeta(t)^2 |e_{yt}(y,t-s)|\Psi(y,s) \,dy\,ds\\
&\le C(E_0+\zeta(t)^2)(1+t)^{-1}.\\
\end{aligned}
$$
Dividing by $(\theta+\psi_1+\psi_2)(x,t)$ and $(1+t)^{-1}$,
respectively, we obtain \eqref{claim} as claimed.
\medbreak

{}From \eqref{claim}, we obtain global existence, with
$\zeta(t)\le 2CE_0$.
{}From the latter bound and the definition of $\zeta$
in \eqref{zeta2} we obtain the first two bounds of \eqref{pointwise}.
It remains to establish the third bound, expressing convergence
of phase $\delta$ to a limiting value $\delta(+\infty)$.

By Lemmas \ref{iniconvolutions}--\ref{convolutions}
together with the previously
obtained bounds \eqref{Nbounds} and $\zeta\le CE_0$,
and the definition \eqref{zeta2} of $\zeta$, the formal limit
$$
 \begin{aligned}
 \delta(+\infty)&:=  \int^\infty_{-\infty}
 e(y,+\infty) U_0(y)\,dy  \\
&\qquad
  +\int^t_0 \int^\infty_{-\infty} |e(y,t-s)(-q,1)^\trans|
  |r(U) (y,s)|\,dy\,ds\\
 &\qquad
 +\int^{+\infty}_0\int^{+\infty}_{-\infty}
 e_y(y,+\infty)(Q(U)+
  \dot \delta U)(y,s)\,dy\,ds \\
&\le
 \int^\infty_{-\infty}
 |e(y,+\infty)|E_0(1+|y|)^{-3/2} \,dy  \\
&\quad +
\int^t_0\int^{+\infty}_{-\infty}|e(y,t-s)(-q,1)^\trans
e^{-\eta y^+}|CE_0\Psi(y,s) \,dy\,ds\\
 &\qquad
 +\int^{+\infty}_0\int^{+\infty}_{-\infty}
 |e_y(y,+\infty)| CE_0 \Psi(y,s) \,dy\,ds \\
&\le CE_0
\end{aligned}
$$
is well-defined, as the sum of absolutely convergent integrals.

Applying Lemmas \ref{iniconvolutions}--\ref{convolutions} a final time,
we obtain
$$
 \begin{aligned}
 |\delta(t)-\delta(+\infty)| &\le \int^\infty_{-\infty}
|e(y,t)-e(y,+\infty)| | U_0(y)| \,dy  \\
&\qquad
  +\int^t_0 \int^\infty_{-\infty}
|\Big(e(y,t-s)- e(y,+\infty)\Big)(-q,1)^\trans|
  |r(U) (y,s)|\,dy\,ds\\
 &\qquad
 +\int^t_0\int^{+\infty}_{-\infty} |e_{y}(y,t-s)-e_y(y,+\infty)|
|(Q(U)+
  \dot \delta U)(y,s)|\,dy\,ds \\
 &\qquad
 +\int_t^{+\infty}\int^{+\infty}_{-\infty} |e_y(y,+\infty)|
|(Q(U)+
  \dot \delta U)(y,s)|\,dy\,ds \\
&\le \int^\infty_{-\infty}
|e(y,t)-e(y,+\infty)| E_0(1+|y|)^{-3/2} \,dy  \\
&\qquad
  +\int^t_0 \int^\infty_{-\infty}
|\Big(e(y,t-s)- e(y,+\infty)\Big)(-q,1)^\trans e^{-\eta y^+}|
CE_0 \Psi(y,s) \,dy\,ds \\
 &\qquad
 +\int^t_0\int^{+\infty}_{-\infty} |e_{y}(y,t-s)-e_y(y,+\infty)|
CE_0 \Psi(y,s) \,dy\,ds \\
 &\qquad
 +\int_t^{+\infty}\int^{+\infty}_{-\infty} |e_y(y,+\infty)|
CE_0 \Psi(y,s)\,dy\,ds \\
&\le CE_0(1+t)^{-1/2},
\end{aligned}
$$
establishing the remaining bound and completing the proof.
\end{proof}

\begin{proof}[Proof of Theorem \ref{nonlin}]
Immediate from Proposition \ref{pwnonlin}, by integration
of bounds \eqref{pointwise}.
\end{proof}

\br\label{zdecay}
\textup{
Proposition \ref{pwnonlin} gives a time-asymptotic description of
perturbation $U$ as a superposition of algebraically decaying
signals propagating along outgoing undamped characteristic directions.
A brief examination reveals that these
consist entirely of {\it fluid dynamical modes}, since reactive modes
propagate always inward from the positive $x$ side,
and as damped outgoing modes on the negative $x$ side.
Recall that fluid modes lie asymptotically along direction
$R_f^\pm=(*,0)^\trans$ with vanishing $z$-component.
Taking account of this fact, together
with the faster decay rate of $(0,1)\tilde G(1,0)^\trans$ stated
in Proposition \ref{greenbounds}, we could by essentially the
same argument used to prove Proposition \ref{pwnonlin}
establish the refined result that {\it the $z$-component of perturbation $U$
decays faster than the $u$-component},
reflecting the physical picture that the fluid is in each case
(weak or strong detonation or deflagration) swept through the traveling wave,
burning completely in the high-temperature region in its wake.
However, we do not determine precise bounds here.
}
\er

\end{document}